\documentclass[notitlepage]{scrartcl}

\usepackage[shortlabels]{enumitem}
\usepackage{amsmath}
\usepackage{amssymb}
\usepackage{amsthm}
\usepackage{abstract}
\usepackage{xcolor}
\usepackage{mathtools}
\usepackage{graphicx}
\usepackage{caption}
\usepackage{subcaption}
\usepackage{float}
\usepackage[margin=1in]{geometry}
\usepackage{hyperref}
\usepackage{cleveref}
\usepackage{tikz}
\usepackage{wrapfig}

\usetikzlibrary{patterns}

\hypersetup{
  colorlinks=true,
  linkcolor=blue,
  filecolor=magenta,
  urlcolor=cyan
}

\newtheorem{theorem}{Theorem}
\newtheorem{thm}{Theorem}[section]

\newtheorem{lemma}[thm]{Lemma}
\newtheorem{proposition}[thm]{Proposition}

\newtheorem{remark}[thm]{Remark}
\newtheorem{question}[thm]{Question}

\newcommand{\whp}{whp}
\newcommand{\defi}{\operatorname{def}}
\newcommand{\odd}{\operatorname{odd}}
\newcommand{\iso}{\operatorname{iso}}

\newcommand{\cO}{\mathcal{O}}
\newcommand{\cE}{\mathcal{E}}

\newcommand{\footremember}[2]{%
  \footnote{#2}%
  \newcounter{#1}%
  \setcounter{#1}{\value{footnote}}%
}

\title{\vspace{-1.5cm}Perfect Matching in Product Graphs\\
and in their Random Subgraphs}
\author{%
Sahar Diskin \footremember{diskin}{\scriptsize School of Mathematical Sciences, Tel Aviv University, Tel Aviv 6997801, Israel. Email: sahardiskin@mail.tau.ac.il.}%
\and
Anna Geisler \footremember{geisler}{\scriptsize Institute of Discrete Mathematics, Graz University of Technology, Steyrergasse 30, 8010 Graz, Austria. Email: geisler@math.tugraz.at.}%
}
\begin{document}

\maketitle
\vspace{-2em}

\begin{abstract}
For $t\in \mathbb{N}$ and ever $i\in [t]$, let $H_i$ be a $d_i$-regular connected graph with $1<|V(H_i)|\le M$ for some integer $M\ge 2$. Let $G=\square_{i=1}^tH_i$ be the $t$-dimensional Cartesian product of $H_1,\ldots, H_t$. We prove that if $t\ge 2\ln M$ then $G$ has a (nearly-)perfect matching. We further show that this bound on the dimension is tight up to a constant factor. 

Then, considering the random graph process on $G$, we generalise the result of Bollob\'as on the binary hypercube $Q^t$, showing that with high probability, the hitting times for minimum degree one, connectivity, and the existence of a (nearly-)perfect matching in the random graph process on $G$ are the same. 
\end{abstract}

\section{Introduction}

\subsection{Background and main results}

Given two graphs $H_1=(V_1,E_1)$ and $H_2=(V_2,E_2)$, the \textit{Cartesian product} $H=H_1\square H_2$ is the graph whose vertex set is $V_1\times V_2$, and for $u_1,v_1\in V_1$ and $u_2,v_2\in V_2$, we have that $\{(u_1,u_2),(v_1,v_2)\}\in E(H)$ either if $u_1=v_1$ and $\{u_2,v_2\}\in E(H_2)$ or if $u_2=v_2$ and $\{u_1,v_1\}\in E(H_1)$. More generally, given $t$ graphs, $H_i,\ldots, H_t$, their Cartesian product $G=\square_{i=1}^tH_i$ is the graph with the vertex set
\begin{align*}
    V\coloneqq \left\{v=(v_1,\ldots, v_t)\colon v_i\in V(H_i)\text{ for all } i\in[t]\right\},
\end{align*}
and the edge set
\begin{align*}
  \left\{uv \colon \begin{array}{l} \text{there is some } i\in [t] \text{ such that }  \{u_i,v_i\}\in E\left(H_i\right)\\
   \qquad\qquad\qquad\quad\qquad \text{and } u_j=v_j \text{ for all } j \neq i \end{array}\right\}.
\end{align*}
We call $H_1,H_2,\ldots, H_t$ the \textit{base graphs of} $G$. Cartesian product graphs arise naturally in many contexts and have received much attention in combinatorics, probability, and computer science. Many classical graphs, which have been extensively studied, are in fact Cartesian product graphs: the $t$-dimensional torus is the Cartesian product of $t$ copies of the cycle $C_k$, the $t$-dimensional grid is the Cartesian product of $t$ copies of the path $P_k$, and the binary $t$-dimensional hypercube $Q^t$ is the Cartesian product of $t$ copies of a single edge $K_2$.
We refer the reader to \cite{HIK11} for a systematic treatment of Cartesian product graphs and related graph products.

A \emph{perfect matching} in a graph $G$ is a collection of vertex disjoint edges which cover the vertex set of $G$. A matching is \emph{nearly-perfect} if it leaves exactly one vertex unmatched. In this paper, we study perfect matchings of high-dimensional product graphs, both in the deterministic and in the random setup. 

We begin with the deterministic part. Observe that if one base graph, say $H_j$, has a perfect matching, then $G=\square_{i=1}^tH_i$ has a perfect matching: for each fixed choice of the other coordinates we may take the corresponding copy of that matching in the $j$-th coordinate.
Kotzig \cite{K79} proved a stronger statement: if all base graphs are regular and either one base graph has a $1$-factorisation, or two base graphs have perfect matchings, then the product has a $1$-factorisation (recall that a $1$-factorisation is a decomposition of the edge set into perfect matchings).
However, when none of the base graphs has a perfect matching, the existence of a perfect matching in the product is not automatic.

Nevertheless, for any two graphs, the proportion of vertices covered by a maximum matching in their product is at least the corresponding proportion in either base graph: one simply takes a maximum matching in that base graph for every fixed choice of the other coordinate. It is therefore natural to ask: can one take sufficiently many products, that is, consider a product graph of high enough dimension, to force a perfect or nearly-perfect matching? If so, what is an optimal bound on the required dimension? Our first result answers the first question in the positive, and gives a \textit{tight} bound on the dimension of a Cartesian product graph that ensures the existence of a perfect matching.

\begin{theorem}\label{th:deterministic}
Let $M\ge2$ and let $t>2 \ln M$ be integers. For every $i\in[t]$, let $H_i$ be a connected $d_i$-regular graph with $1<|V(H_i)|\le M$, and write $G=\square_{i=1}^t H_i$. Then $G$ has a (nearly-)perfect matching.

Furthermore, there is a sequence of such graphs $(H_i)_{i \in [t]}$ such that for $t<\ln M/3$ the product $G=\square_{i=1}^t H_i$ does not have a perfect matching.
\end{theorem}


In the following remark, we show that some regularity assumption is necessary for a perfect matching to appear. 

\begin{remark}\label{rem:stars}
Let $G=\square_{i=1}^tK_{1,s}$.
The graph $G$ is bipartite; writing $\cE$ and $\cO$ for the sides according to the parity of the number of leaf coordinates, we have
\[
  |\cE|=\sum_{\substack{0\le k\le t\\k\text{ even}}}\binom{t}{k}s^k,
  \qquad
  |\cO|=\sum_{\substack{0\le k\le t\\k\text{ odd}}}\binom{t}{k}s^k.
\]
Consequently,
\[
  \bigl||\cE|-|\cO|\bigr|=|1-s|^t.
\]
For $s\ge3$ this imbalance is larger than one, so $G$ has neither a perfect nor a nearly-perfect matching.
\end{remark}

We now turn to the random part.
Given a graph $G$ and probability $p$, we form the percolated random subgraph $G_p\subseteq G$ by including every edge of $G$ independently with probability $p$ (note that $G(n,p)$ is then $(K_n)_p$). 

The study of $Q^t_p$ was initiated by Sapo\v{z}enko \cite{S67} and by Burtin \cite{B77}, who showed that the (sharp) threshold for connectivity is $p^*=\frac{1}{2}$: when $p<\frac{1}{2}$ is a constant, with high probability (whp), that is, with probability tending to one as $t$ tends to infinity, $Q^t_p$ is disconnected, whereas for $p>\frac{1}{2}$ constant, whp $Q^t_p$ is connected. Erd\H{o}s and Spencer \cite{ES79} conjectured that $Q^t_p$ undergoes a phase transition with respect to its component structure, that is, the typical emergence of a giant component (a connected component containing a linear fraction of the vertices) around $p=\frac{1}{t}$, similar to that of $G(n,p)$ around $p=\frac{1}{n}$. This conjecture was confirmed by Ajtai, Koml\'os, and Szemer\'edi \cite{AKS81}, with subsequent work by Bollob\'as, Kohayakawa, and \L{}uczak~\cite{BKL92}. We refer the reader to \cite{K23} for a modern short proof of this result.

In recent years, there has been an effort to generalise these results to a wider family of product graphs. Lichev \cite{L22} gave sufficient conditions, in terms of the base graphs, for the typical emergence of a giant component in bond percolation on high-dimensional product graphs. Diskin, Erde, Kang, and Krivelevich \cite{DEKK23} improved upon this, giving sufficient and tight conditions for the typical emergence of a giant component. Furthermore, they showed that assuming the base graphs are regular, one can give a rather precise description, similar to that in $G(n,p)$, both of the typical component structure~\cite{DEKK22}, and of the asymptotic combinatorial properties of the giant component~\cite{DEKK24}. 

Returning to the dense regime, that is, when $p$ is constant, connectivity and the existence of a perfect matching in $Q^t_p$ have been studied in detail. In particular, Bollob\'as \cite{B90} obtained a \textit{hitting time} result for the \textit{random graph process} on $Q^t$. Given a graph $\Gamma$, the random graph process on $\Gamma$ is defined as a random sequence of nested graphs $\Gamma(0)\subseteq \ldots \subseteq \Gamma(|E(\Gamma)|)$ together with an ordering $\sigma$ on $E(\Gamma)$, chosen uniformly at random from among all $|E(\Gamma)|!$ such orderings. We set $\Gamma(0)$ to be the empty graph on $V(\Gamma)$. Given $\Gamma(i)$, with $0\le i < |E(\Gamma)|$, we form $\Gamma(i+1)$ by adding the $(i+1)$-st edge, according to the ordering $\sigma$, to $\Gamma(i)$. The hitting time of a monotone increasing, non-empty graph property $\mathcal{P}$, is the random variable equal to the index $\tau$ for which $\Gamma(\tau)\in \mathcal{P}$, but $\Gamma(\tau-1)\notin\mathcal{P}$. Note that having minimum degree one, connectivity, and the existence of a perfect matching are all monotone increasing properties. Furthermore, observe that for a graph $\Gamma$ to be connected or to contain a perfect matching, the minimum degree of $\Gamma$ has to be at least one. A classical result of Erd\H{o}s and R\'enyi \cite{ER66}, and of Bollob\'as and Thomason \cite{BT85}, is that in the random process on $K_n$, whp the hitting time for minimum degree one, connectivity, and the existence of a perfect matching are the same. In 1990, Bollob\'as \cite{B90} showed that the same phenomenon holds in the random graph process on $Q^t$.
Subsequent work by Joos determined the threshold for connectivity for Cartesian powers of graphs \cite{J12, J15}, that is, for $G=\square_{i=1}^tH$ where all the base graphs $H_i=H$ are the same.


The corresponding statement for Cartesian products whose base graphs have bounded order is our main probabilistic result. Here and below, an event holds \emph{with high probability}, abbreviated \whp, if its probability tends to one as $t\to\infty$ with $M$ fixed.

\begin{theorem}\label{th:hitting-time}
Let $M\ge2$ be fixed, and for every $i\in[t]$ let $H_i$ be a connected $d_i$-regular graph with $1<|V(H_i)|\le M$.
Let $G=\square_{i=1}^tH_i$, and suppose $|V(G)|$ is even.
In the random graph process on $G$, let $\tau_1$, $\tau_2$, and $\tau_3$ be the hitting times for minimum degree one, connectivity, and the existence of a perfect matching, respectively.
Then, \whp, $  \tau_1=\tau_2=\tau_3$.
\end{theorem}

We use a standard coupling of the random graph process with the random subgraph formed by percolation: assign to every $e\in E(G)$ a label $\xi_e$ that is uniform on $[0,1]$, and set
\[
  G_p=\bigl(V(G),\{e\in E(G):\xi_e\le p\}\bigr).
\]
The labels are almost surely distinct, and $(G_p)_{0\le p\le1}$ then naturally couples with the random graph process.
%
For a graph $\Gamma$, write $\iso(\Gamma)$ for the number of its isolated vertices.
We work on the following interval around the disappearance of the last isolated vertices.
Write $n=|V(G)|$ and $d=d(G)$, let $\omega=\omega(d)$ satisfy
\[
  \omega\to \infty
  \qquad\text{and}\qquad
  \omega\le d^{1/100},
\]
and define $p_-$ and $p_+$ by
\begin{equation}\label{eq:window}
  (1-p_-)^d=\frac{\omega}{n},
  \qquad
  (1-p_+)^d=\frac{1}{\omega n}.
\end{equation}
We are now ready to state the following structural result for random subgraphs of $G$.
\begin{theorem}\label{th:window-connectivity}
Let $M\ge2$ be fixed, and for every $i\in[t]$ let $H_i$ be a connected $d_i$-regular graph with $1<|V(H_i)|\le M$.
Let $G=\square_{i=1}^tH_i$, and suppose $|V(G)|$ is even. Let $p \in [p_-, p_+]$.  Then, \whp, the following statements hold simultaneously.
\begin{enumerate}[(a)]
  \item In $G_{p_-}$ there is a unique component containing every non-isolated vertex. The isolated vertices are pairwise at distance at least three in $G$, and
  \[
    1\le \iso(G_{p_-})\le \sum_{i=1}^td_i.
  \]
  \item The graph $G_{p_+}$ has no isolated vertices.
  \item For any $p\in [p_-, p_+]$, if $G_p$ has minimum degree at least one, then it has a perfect matching. 
\end{enumerate}
\end{theorem}

Note, that \Cref{th:window-connectivity} readily implies \Cref{th:hitting-time}: Part (a) and (b) of \Cref{th:window-connectivity} imply that with high probability $\tau_1\in[p_-,p_+]$ and $G_{\tau_1}$ is connected, so $\tau_2=\tau_1$.
Since $\delta(G_{\tau_1})\ge1$, part (c) of \Cref{th:window-connectivity} gives a perfect matching in $G_{\tau_1}$.
Having a perfect matching is an increasing property and implies minimum degree at least one; therefore $\tau_3=\tau_1$ as well.

One of the key challenges in general product graphs, as compared to the hypercube, is the lack of bipartitedness. In particular, the existence of a perfect matching cannot be proved by verifying Hall's condition; instead, we use the Tutte--Berge formula and the Edmonds-Gallai decomposition. 

\paragraph{Structure of the paper}

The paper is structured as follows. In Section \ref{s:prelim}, we collect some notation and results that will be of use throughout the paper. In \Cref{s:outline} we give a proof outline containing the key ideas. In \Cref{sec:deterministic-matching} we cover the deterministic setting and prove \Cref{th:deterministic}. Switching to the random setting, we treat connectivity in \Cref{s:connectivity}, and the existence of a perfect matching in \Cref{s:matching}. Finally, we conclude with a discussion and some open problems in \Cref{s:discussion}.

\section{Preliminaries}\label{s:prelim}

\paragraph{Notation and definitions.}
Throughout, $M\ge2$ is fixed, $(H_i)_{i=1}^t$ is a sequence of connected $d_i$-regular graphs with $1<|V(H_i)|\le M$, and
\[
  G=(V,E)=\square_{i=1}^tH_i,
  \qquad n=|V|,
  \qquad d=d(G)=\sum_{i=1}^td_i.
\]
We call $t$ the \emph{dimension} of $G$, and, for $u=(u_1,\ldots,u_t)\in V$, we call $u_i$ the $i$-th coordinate of $u$.
Since $2^t\le n\le M^t$ and $t\le d\le(M-1)t$, we have $d=\Theta_M(\ln n)$. More precisely, we will use repeatedly that $\ln n \geq \frac{\ln 2}{M-1} d$.
For a graph $\Gamma$ and a set $A\subseteq V(\Gamma)$, write $E_\Gamma(A)$ for the set of edges of $\Gamma$ with both endpoints in $A$, and set $e_\Gamma(A)\coloneqq |E_\Gamma(A)|$.
For disjoint sets $A,B\subseteq V(\Gamma)$, write $E_\Gamma(A,B)$ for the set of edges of $\Gamma$ between $A$ and $B$, and set $e_\Gamma(A,B)\coloneqq |E_\Gamma(A,B)|$.
We abbreviate $A^c=V(\Gamma)\setminus A$ when the graph in question is clear. Furthermore, we denote by $\lambda(\Gamma)$ the minimum number of edges whose removal disconnects $\Gamma$ and call this the edge-connectivity of $\Gamma$. The deficit of $\Gamma$, abbreviated by $\defi(\Gamma)$, is the number of vertices left unmatched by a maximum matching.
Finally, we write $\odd(\Gamma)$ for the number of odd-order components in a graph $\Gamma$.

\paragraph{External results}
We use the following isoperimetric inequalities for product graphs with bounded-order base graphs; they follow from Theorems~1 and~2 of \cite{DEKK24} and Section~3.4 of \cite{DS24}.

\begin{thm}\label{th:expansion}
Let $M\ge2$ be fixed, and for every $i\in[t]$ let $H_i$ be a connected $d_i$-regular graph with $1<|V(H_i)|\le M$.
Let $G=(V, E)=\square_{i=1}^tH_i$. For every non-empty $A\subseteq V$,
\begin{align}
  e_G(A,A^c)&\ge \frac{e}{M}|A|\ln\frac{n}{|A|},\label{eq:expansion-log}\\
  e_G(A,A^c)&\ge |A|\bigl(d-(M-1)\log_M|A|\bigr).\label{eq:expansion-local}
\end{align}
\end{thm}

We also use the standard rooted-tree bound (see, for example, \cite[Lemma~2]{BFM98}).

\begin{lemma}\label{lem:rooted-trees}
Let $\Gamma$ be a graph of maximum degree at most $d$, let $k\ge1$ be an integer, and fix some $v\in V(\Gamma)$.
The number of $k$-vertex trees in $\Gamma$ containing $v$ is at most $(ed)^{k-1}$.
\end{lemma}

\paragraph{The Tutte-Berge formula}

The Tutte-Berge formula \cite{B62} is our first tool in order to treat the existence of perfect matchings. It implies that the number of unmatched vertices in a maximum matching in a graph $\Gamma$ is equal to 
\begin{align} \label{eq:Tutte-Berge}
    \defi(\Gamma)= \max_{U\subseteq V(\Gamma)} \left\{\odd\left(\Gamma[V(\Gamma)\setminus U]\right)- |U|\right\},
\end{align}
where we recall that $\odd\left(\Gamma[V(\Gamma)\setminus U]\right)$ is the number of connected components with an odd number of vertices in $\Gamma[V(\Gamma)\setminus U]$. In particular, a graph $\Gamma$ has a perfect matching if and only if for every subset $U\subseteq V(\Gamma)$, the subgraph $\Gamma[V(\Gamma)\setminus U]$ has at most $|U|$ connected components with an odd number of vertices. Note here that if $\Gamma$ has isolated vertices, then choosing $U=\varnothing$ witnesses that there is no perfect matching in $\Gamma$.

\paragraph{The Gallai-Edmonds decomposition} \label{s:Gallai}
The second tool provides us with a structural description of any maximum matching in $\Gamma$. We define
\begin{align*}
 D_\Gamma&\coloneqq\{v\in V(\Gamma):\text{some maximum matching of $\Gamma$ leaves $v$ unmatched}\},\\
 A_\Gamma&\coloneqq N_\Gamma(D_\Gamma)\setminus D_\Gamma,\\
 C_\Gamma&\coloneqq V(\Gamma)\setminus(D_\Gamma\cup A_\Gamma).
\end{align*}
The uniquely determined partition
\[
 V(\Gamma)=D_\Gamma\sqcup A_\Gamma\sqcup C_\Gamma
\]
is the \emph{Gallai--Edmonds decomposition} of $\Gamma$.
A graph is called \emph{factor-critical} if deleting any one of its vertices leaves a graph with a perfect matching.
Writing $m_{\gamma}$ for the number of components of $\Gamma[D_\Gamma]$ and enumerating them by $D_1,\ldots,D_{m_\Gamma}$, the Gallai--Edmonds theorem \cite{edmonds1965paths,gallai1964maximale} states that
\begin{enumerate}[(i)]
 \item every $D_i$ is factor-critical;
 \item $\Gamma[C_\Gamma]$ has a perfect matching;
 \item every maximum matching of $\Gamma$ consists of a nearly-perfect matching in each $D_i$, a perfect matching of $\Gamma[C_\Gamma]$, and $|A_\Gamma|$ edges that match the vertices of $A_\Gamma$ to pairwise distinct components $D_i$.
\end{enumerate}
In particular,
\begin{align*} \label{eq:Gallai-Edmonds}
     \defi(\Gamma)=m_\Gamma-|A_\Gamma|.
\end{align*}

Moreover, there is no edge of $\Gamma$ between $D_\Gamma$ and $C_\Gamma$.
We refer to \cite[Chapter~3]{LP86} for the theorem and further background.


\section{Proof outline}\label{s:outline}

\paragraph{Deterministic setting: Theorem~\ref{th:deterministic}.}
If one base graph has a perfect matching, then automatically the product has a perfect matching. However, even if none of the base graphs has a perfect matching, the proportion of vertices covered by a maximal matching is preserved and we hope to increase it when taking products of higher and higher dimensions. We prove this by using that the edge-connectivity of a Cartesian product graph indeed increases with its dimension. With this at hand, let $U \subseteq V$ be an obstruction to a perfect matching in the sense of the Tutte-Berge formula, that is $G - U$ contains more than $|U|$ odd components. The edge-connectivity from \Cref{lem:product-connectivity} implies that the edge-boundary of $U$ is of size at least $d-1$, and a parity argument shows then shows that in fact this is of size least $d$. All these boundary edges must have one vertex in $U$ and therefore, summing over the odd components in $G-U$, there can be at most $|U|$ of them since the graph under consideration is $d$-regular. By the Tutte-Berge formula, we conclude that such a product of sufficiently high dimension contains a perfect matching deterministically.

\paragraph{Random setting: Connectivity}

We prove a statement that holds simultaneously for every $G_p$ with $p \in [p_-,p+]$. In order to treat the random graph $G_p$ we use a classical sprinkling argument, that is, we expose the random edges in two independent rounds. The second, much sparser round, is referred to as the sprinkling step. After the first round, first-moment calculations yield the following properties of the resulting graph whp: its isolated vertices are at distance at least three in $G$ and it has no component whose order lies in between $2$ and $d^{20}$, see \Cref{lem:isolates-far,lem:component-gap}. It remains to merge the components of order at least $d^{20}$. For every partition of these components into two non-empty parts, the expansion properties of $G$ yield many short edge-disjoint paths from one part to another. The probability that the second round contains none of these paths is very small, and we conclude using a union bound argument. This proves that whp in the range of $p \in [p_-, p+]$ there is a unique giant component and all vertices outside are isolated. Combining this with the fact that these vertices are at distance at least three in $G$, this implies the hitting time result for connectivity and minimum degree one. Indeed, each such isolated vertex joins this giant component once it receives an edge, and the isolated vertices disappear within that range of parameters whp. Thus, the graph becomes connected exactly at the same time the last isolated vertex receives an edge whp.

\paragraph{Random setting: perfect matching}

Again, we prove a statement that holds simultaneously for every $G_p$ with $p \in [p_-,p+]$ under the additional assumption that $\delta(G_p) \geq 1$. By the first part, we my assume that such a random graph is connected whp. Our claim is that at this point it also contains a perfect matching. Suppose otherwise and consider the Edmond-Gallai decomposition $(D, A, C)$ of $G_p$. If $D_1, ..., D_m$ are the components of $G_p[D]$ and writing $a\coloneqq |A|$, then
\[m-a = def(G_p) \geq 2\]
since $G$ has even order, so that every maximum matching leaves at least two vertices unmatched. The components $D_i$ are odd, every component of $G_p[C]$ is even and there is no edge between $D$ and $C$. Consequently,
\[odd(G_p-A)=m \geq a+2\]
and $A$ is an obstruction to a perfect matching in the sense of the Tutte-Berge formula. Our main task is thus to prove the following two statements:

\begin{enumerate}
    \item Whp, there are no small sets $U$ forming such an obstruction.
\item Whp, the set $A$ from the Edmonds-Gallai decomposition of $G_p$ is small.
\end{enumerate}

For the first statement, we first prove any such minimal set $U$ has neighbours in at least three odd components of $G_p-U$, see \Cref{lem:three_edges}. Then we crucially use the expansion properties of the graph to show that a large number of edges must be absent in $G_p$ to allow for such an obstruction $U$. This is combined with a union bound over the number of such obstructions depending on its size and expansion. Crucially, we prove an improved bound on the number of sets with small edge-boundary, see \Cref{lem:small-boundary-sets}. This bound makes an explicit use of the Cartesian product structure: indeed, sets with bad expansion leave a `fingerprint' on many coordinates, limiting their number.

It remains to study the Edmonds-Gallai decomposition of $G_p$ to conclude that whp the part $A$ is not too large. Again, note that every edge joining different components of $G_p [D]$ and every edge from $D$ to $C$ must be absent in $G_p$. A degree count expresses the number of these absent edges in terms of the deficiency $def(G_p)$, the edges inside $A$, the boundary of $C$ and the edge boundaries of components of $G_p[D]$, see \eqref{eq:combined-LB-h}. Using expansion properties of the graph, we lower-bound the number of edges leaving these components.

We then distinguish between the case when $G$ is bipartite and non-bipartite. If $G$ is bipartite, then every factor-critical bipartite graph consists of a single vertex, so all components of $G_p[D]$ are singletons. Using that the bipartition is balanced, we find a specific set whose internal edges must be absent and counting via choosing rooted spanning trees allows for a union bound.

If $G$ is non-bipartite, some component of $G_p[D]$ could contain more than one vertex. Thus, we distinguish between the singleton and non-singleton components of  $D$. After these components are fixed, it remains to choose the vertices in $A$ and the set $C$. 
It remains crucial to balance a lower bound on the number of absent edges with an upper bound on the number of choices for such a structure. In particular, in \Cref{lem:two-colour-count}, we bound the number of choices for a bipartition of a large vertex set with few internal edges in $G_p$.

Recalling the coupling of the random graph process with $G_p$, we note that all bounds require only that absent edges have labels larger than $p_-$ and present edges have labels smaller than $p^+$. Therefore, they hold simultaneously throughout the interval $p \in [p_-, p^+]$, which contains the hitting time for minimum degree one (and connectivity) whp. The fact that at this point $G_p$ also has a perfect matching shows the equality of the three hitting times. 

\section{Deterministic setting: Existence of a perfect matching} \label{sec:deterministic-matching}

We will deduce the existence statement in \Cref{th:deterministic} from edge-connectivity properties of the Cartesian product.

\begin{lemma}\label{lem:product-connectivity}
Let $G=\square_{i=1}^t H_i$, where every $H_i$ is connected and $d_i$-regular.
Let $M_i=|V(H_i)|$ and $d=\sum_{i=1}^t d_i$.
Then
\[
  \lambda(G)=
  \min\left\{d,\ \min_{i\in[t]}\lambda(H_i)\prod_{j\ne i}M_j\right\}.
\]
\end{lemma}

\begin{proof}
The displayed identity is the standard edge-connectivity formula for Cartesian products; see \cite[Chapter~5]{HIK11}.
\end{proof}

We will use the following lemma, which follows almost immediately from the Tutte--Berge formula; we include its proof for the sake of completeness.
\begin{lemma}\label{lem:regular-matching}
Let $d \geq 1$ be an integer and let $\Gamma$ be a connected $d$-regular graph. Suppose that $\lambda(\Gamma)\ge d-1$.
Then $\Gamma$ has a perfect matching when $|V(\Gamma)|$ is even and a nearly-perfect matching when $|V(\Gamma)|$ is odd.
\end{lemma}

\begin{proof}
For $U\subseteq V(\Gamma)$, let $Q_1,\ldots,Q_q$ be the odd components of $\Gamma-U$ such that $q=\odd(\Gamma - U)$.
If $U=\varnothing$, then $q=0$ when $|V(\Gamma)|$ is even and $q=1$ when it is odd.
Suppose now that $U\ne\varnothing$.
Every $Q_i$ is a proper non-empty set, so by the assumption that $\lambda(\Gamma)\geq d-1$,
\[
  e_\Gamma(Q_i,Q_i^c)\ge d-1.
\]
Since $|Q_i|$ is odd and $\Gamma$ is $d$-regular (implying that $d|Q_i|= 2 e_{\Gamma}(Q_i) + e_{\Gamma}(Q_i, Q_i^c)$), we have
\[
  e_\Gamma(Q_i,Q_i^c)\equiv d|Q_i|\equiv d\pmod 2.
\]
The cut has the same parity as $d$, whereas $d-1$ has the opposite parity; hence
\[
  e_\Gamma(Q_i,U)=e_\Gamma(Q_i,Q_i^c)\ge d.
\]
Summing over the odd components yields
\[
  qd\le\sum_{i=1}^q e_\Gamma(Q_i,U)\le d|U|,
\]
and therefore $q\le|U|$.
The Tutte--Berge formula \eqref{eq:Tutte-Berge} implies that all vertices are covered by a maximal matching if $|V(\Gamma)|$ is even and all but at most one vertex are covered by a maximal matching if $|V(\Gamma)|$ is odd.
\end{proof}

Combining these two lemmas yields a perfect matching in product graphs of sufficiently high dimension.

\begin{proof}[Proof of \Cref{th:deterministic}]
Let $G=\square_{i=1}^t H_i$. If a base graph $H_i$ has a perfect matching $M_i$, then for each fixed choice of all coordinates other than the $i$-th, take the corresponding copy of $M_i$. The resulting matchings are vertex-disjoint and cover $V(G)$.

Otherwise, any base graph $H_i$ has at least three vertices. Thus, by \Cref{lem:product-connectivity} we obtain that
\begin{align} \label{eq:edge-conn}
      \lambda(G)
  =\min\left\{d,\ \min_{i\in[t]}\lambda(H_i)\prod_{j\ne i}|V(H_j)|\right\} \geq \min\{d, 3^{t-1}\}.
\end{align}
Using that $d \leq t(M-1)$ and $t > 2\ln M$ we obtain that \eqref{eq:edge-conn} is at least $d-1$.
Applying \Cref{lem:regular-matching}, we obtain a (nearly-) perfect matching. This concludes the first part of the theorem.

For the second part, we provide an explicit construction where for $t<\ln M/3$ there is no perfect matching. Indeed, let $M$ be sufficiently large, and let $m$ be the largest odd integer with $(m+1)^2\le M$. 
Construct a graph $H'$ as follows. Consider $K_{m-1}$ and three distinct vertices $x,y$, and $z$. Remove one perfect matching from $K_{m-1}$, draw the edge $yz$, and draw all edges between $x,y,z$ and $K_{m-1}$. Thus the vertex $x$ has degree $m-1$, and every other vertex in $H'$ has degree $m$. We then construct $H$ by taking $m$ copies of $H'$ and a distinct vertex $\tilde x$, and drawing the edges between $\tilde x$ and the copy of $x$ in each copy of $H'$, see \Cref{fig:construction-tight}. Thus $H$ is an $m$-regular graph on $(m+1)^2$ vertices. Moreover, since $m$ is odd, $H\setminus \{\tilde x\}$ has $m$ odd components, and thus $H$ has no perfect matching by Tutte's Theorem. Finally, construct $G$ by taking the Cartesian product of $H$ with $(t-1)$ triangles $K_3$. The base graphs all have size at most $M$, are connected and regular. Consider the set $U=\{\tilde x\}\times V(\square_{i=1}^{t-1}K_3)$. $G-U$ has exactly $m$ odd components, so if $|U|=3^{t-1}<m\le \sqrt{M}$, $G$ has no perfect matching. In particular, any $t<\ln M/3$ suffices.
\end{proof}

\begin{figure}
    \centering

\tikzset{every picture/.style={line width=0.75pt}} 

\begin{tikzpicture}[x=0.75pt,y=0.75pt,yscale=-1,xscale=1]

\draw   (240,35.55) .. controls (240,27.51) and (246.51,21) .. (254.55,21) -- (342.35,21) .. controls (350.39,21) and (356.9,27.51) .. (356.9,35.55) -- (356.9,79.21) .. controls (356.9,87.24) and (350.39,93.76) .. (342.35,93.76) -- (254.55,93.76) .. controls (246.51,93.76) and (240,87.24) .. (240,79.21) -- cycle ;
\draw  [fill={rgb, 255:red, 0; green, 0; blue, 0 }  ,fill opacity=1 ] (295.32,128.27) .. controls (295.32,126.21) and (297.09,124.54) .. (299.28,124.54) .. controls (301.46,124.54) and (303.23,126.21) .. (303.23,128.27) .. controls (303.23,130.33) and (301.46,132) .. (299.28,132) .. controls (297.09,132) and (295.32,130.33) .. (295.32,128.27) -- cycle ;
\draw  [fill={rgb, 255:red, 0; green, 0; blue, 0 }  ,fill opacity=1 ] (367.44,71.37) .. controls (367.44,69.31) and (369.21,67.64) .. (371.39,67.64) .. controls (373.58,67.64) and (375.35,69.31) .. (375.35,71.37) .. controls (375.35,73.43) and (373.58,75.1) .. (371.39,75.1) .. controls (369.21,75.1) and (367.44,73.43) .. (367.44,71.37) -- cycle ;
\draw  [fill={rgb, 255:red, 0; green, 0; blue, 0 }  ,fill opacity=1 ] (368.43,34.99) .. controls (368.43,32.93) and (370.2,31.26) .. (372.38,31.26) .. controls (374.56,31.26) and (376.33,32.93) .. (376.33,34.99) .. controls (376.33,37.05) and (374.56,38.72) .. (372.38,38.72) .. controls (370.2,38.72) and (368.43,37.05) .. (368.43,34.99) -- cycle ;
\draw [line width=1.5]    (372.38,34.99) -- (371.39,71.37) ;
\draw [color={rgb, 255:red, 128; green, 128; blue, 128 }  ,draw opacity=1 ]   (254.16,86.29) -- (299.28,128.27) ;
\draw [color={rgb, 255:red, 128; green, 128; blue, 128 }  ,draw opacity=1 ]   (299.28,128.27) -- (346.04,86.29) ;
\draw [color={rgb, 255:red, 128; green, 128; blue, 128 }  ,draw opacity=1 ]   (372.38,34.99) -- (331.33,28) ;
\draw [color={rgb, 255:red, 128; green, 128; blue, 128 }  ,draw opacity=1 ]   (371.39,71.37) -- (334.33,85) ;
\draw [color={rgb, 255:red, 128; green, 128; blue, 128 }  ,draw opacity=1 ]   (372.38,34.99) -- (335.33,52) ;
\draw [color={rgb, 255:red, 128; green, 128; blue, 128 }  ,draw opacity=1 ]   (371.39,71.37) -- (337.33,60) ;
\draw [line width=1.5]    (299.28,128.27) -- (346,188) ;
\draw  [fill={rgb, 255:red, 0; green, 0; blue, 0 }  ,fill opacity=1 ] (342,188) .. controls (342,185.79) and (343.79,184) .. (346,184) .. controls (348.21,184) and (350,185.79) .. (350,188) .. controls (350,190.21) and (348.21,192) .. (346,192) .. controls (343.79,192) and (342,190.21) .. (342,188) -- cycle ;
\draw   (80,37.55) .. controls (80,29.51) and (86.51,23) .. (94.55,23) -- (182.35,23) .. controls (190.39,23) and (196.9,29.51) .. (196.9,37.55) -- (196.9,81.21) .. controls (196.9,89.24) and (190.39,95.76) .. (182.35,95.76) -- (94.55,95.76) .. controls (86.51,95.76) and (80,89.24) .. (80,81.21) -- cycle ;
\draw  [fill={rgb, 255:red, 0; green, 0; blue, 0 }  ,fill opacity=1 ] (135.32,130.27) .. controls (135.32,128.21) and (137.09,126.54) .. (139.28,126.54) .. controls (141.46,126.54) and (143.23,128.21) .. (143.23,130.27) .. controls (143.23,132.33) and (141.46,134) .. (139.28,134) .. controls (137.09,134) and (135.32,132.33) .. (135.32,130.27) -- cycle ;
\draw  [fill={rgb, 255:red, 0; green, 0; blue, 0 }  ,fill opacity=1 ] (207.44,73.37) .. controls (207.44,71.31) and (209.21,69.64) .. (211.39,69.64) .. controls (213.58,69.64) and (215.35,71.31) .. (215.35,73.37) .. controls (215.35,75.43) and (213.58,77.1) .. (211.39,77.1) .. controls (209.21,77.1) and (207.44,75.43) .. (207.44,73.37) -- cycle ;
\draw  [fill={rgb, 255:red, 0; green, 0; blue, 0 }  ,fill opacity=1 ] (208.43,36.99) .. controls (208.43,34.93) and (210.2,33.26) .. (212.38,33.26) .. controls (214.56,33.26) and (216.33,34.93) .. (216.33,36.99) .. controls (216.33,39.05) and (214.56,40.72) .. (212.38,40.72) .. controls (210.2,40.72) and (208.43,39.05) .. (208.43,36.99) -- cycle ;
\draw [line width=1.5]    (212.38,36.99) -- (211.39,73.37) ;
\draw [color={rgb, 255:red, 128; green, 128; blue, 128 }  ,draw opacity=1 ]   (94.16,88.29) -- (139.28,130.27) ;
\draw [color={rgb, 255:red, 128; green, 128; blue, 128 }  ,draw opacity=1 ]   (139.28,130.27) -- (186.04,88.29) ;
\draw [color={rgb, 255:red, 128; green, 128; blue, 128 }  ,draw opacity=1 ]   (212.38,36.99) -- (171.33,30) ;
\draw [color={rgb, 255:red, 128; green, 128; blue, 128 }  ,draw opacity=1 ]   (211.39,73.37) -- (174.33,87) ;
\draw [color={rgb, 255:red, 128; green, 128; blue, 128 }  ,draw opacity=1 ]   (212.38,36.99) -- (175.33,54) ;
\draw [color={rgb, 255:red, 128; green, 128; blue, 128 }  ,draw opacity=1 ]   (211.39,73.37) -- (177.33,62) ;
\draw [line width=1.5]    (139.28,130.27) -- (346,188) ;
\draw [color={rgb, 255:red, 0; green, 0; blue, 0 }  ,draw opacity=1 ][line width=1.5]    (346,188) -- (444.33,126) ;

\draw (363,188) node [anchor=north west][inner sep=0.75pt]   [align=left] {$\displaystyle \tilde{x}$};
\draw (421,57) node [anchor=north west][inner sep=0.75pt]   [align=left] {...};
\draw (312,116) node [anchor=north west][inner sep=0.75pt]   [align=left] {$\displaystyle x$};
\draw (381,63) node [anchor=north west][inner sep=0.75pt]   [align=left] {$\displaystyle z$};
\draw (382,21) node [anchor=north west][inner sep=0.75pt]   [align=left] {$\displaystyle y$};
\draw (250,47) node [anchor=north west][inner sep=0.75pt]   [align=left] {$\displaystyle K_{m-1} \setminus PM$};
\draw (152,118) node [anchor=north west][inner sep=0.75pt]   [align=left] {$\displaystyle x$};
\draw (221,65) node [anchor=north west][inner sep=0.75pt]   [align=left] {$\displaystyle z$};
\draw (222,23) node [anchor=north west][inner sep=0.75pt]   [align=left] {$\displaystyle y$};
\draw (90,49) node [anchor=north west][inner sep=0.75pt]   [align=left] {$\displaystyle K_{m-1} \setminus PM$};

\end{tikzpicture}

    \caption{The construction of $H$.}
    \label{fig:construction-tight}
\end{figure}

\section{Connectivity in the random graph process}\label{s:connectivity}

Throughout this section, let $p$ satisfy
\begin{equation}\label{eq:mu-range}
  d^{-1/50}\le n(1-p)^d\le d^{1/50}.
\end{equation}
We will expose the random edges in two independent rounds.
In the first round every edge is present with probability $p_1$, and in the second round every edge is present with probability $p_2=d^{-2}$, where $(1-p_1)(1-p_2)=1-p$.

Note that \eqref{eq:mu-range} includes $p=p_-$ and $p=p_+$, as well as $p_1$ from the two-round exposure.
In view of $d=\Theta_M(\ln n)$, \eqref{eq:mu-range} implies that both $p$ and $1-p$ are bounded below by positive constants depending only on $M$.
We begin by showing that isolated vertices are typically far apart in $G$.

\begin{lemma}\label{lem:isolates-far}
Whp, every two isolated vertices of $G_p$ are at distance at least three in $G$.
\end{lemma}

\begin{proof}
There are at most $nd$ ordered adjacent pairs and at most $nd^2$ ordered pairs at distance two.
If distinct vertices $u,v$ are not adjacent, the sets of edges incident with them are disjoint, and hence
\[
  \mathbb{P}(u,v\text{ are isolated in }G_p)=(1-p)^{2d}.
\]
If $u$ and $v$ are adjacent, their incident edge sets overlap in one edge, so this probability is
\[
  (1-p)^{2d-1}\le \alpha_M(1-p)^{2d},
\]
for some constant $\alpha_M$.
The expected number of ordered pairs of isolated vertices at distance at most two is therefore at most
\[
  \alpha_M nd^2(1-p)^{2d}   \le \frac{\alpha_M d^{2+1/25}}{n}=o(1),
\]
where we used that $p$ satisfies \eqref{eq:mu-range}.
We conclude with Markov's inequality.
\end{proof}

Let us further show the following gap statement, which follows from the expansion estimates in \Cref{th:expansion} and the rooted-tree bound in \Cref{lem:rooted-trees}.

\begin{lemma}\label{lem:component-gap}
Whp, $G_p$ has no component whose order lies in $[2,d^{20}]$.
\end{lemma}

\begin{proof}
Fix $2\le k\le d^{20}$.
For every set $K\subseteq V$ of size $k$, the expansion from \eqref{eq:expansion-local} gives
\[
  e_G(K,K^c)\ge k\bigl(d-(M-1)\log_Mk\bigr)\ge0.9kd.
\]
If $K$ is a component of $G_p$, it contains a spanning tree all of whose edges are present, while every edge of $E_G(K,K^c)$ is absent.
Thus, using \Cref{lem:rooted-trees}, we obtain
\begin{align*}
  \mathbb{P}(G_p\text{ has a component of order }k)
  \le n(ed)^{k-1}(1-p)^{0.9kd} \le n^{1-0.8k},
\end{align*}
Where the last inequality follows since $2\le k\le d^{20}$,  $n (1-p)^d\le d^{1/50}$ and $d=\Theta_M(\ln n)$.
Summing over $k$ shows that whp there is no component whose order lies in $[2, d^{20}]$.
\end{proof}

We are now ready to prove the key connectivity statement.

\begin{proposition}\label{prop:unique-component}
Whp, $G_p$ has a unique component containing every non-isolated vertex.
Its isolated vertices are pairwise at distance at least three in $G$.
\end{proposition}

\begin{proof}
Let $p_2=d^{-2}$, and choose $p_1$ so that
\[
  (1-p_1)(1-p_2)=1-p.
\]
The union of independent copies of $G_{p_1}$ and $G_{p_2}$ has the distribution of $G_p$, and
\[
  n(1-p_1)^d=\frac{n(1-p)^d}{(1-d^{-2})^d}=(1+o(1))n (1-p)^d.
\]
The proofs of \Cref{lem:isolates-far,lem:component-gap} hold in the extended regime given in \eqref{eq:mu-range} and thus they apply to $G_{p_1}$.
Also, the expected number of isolated vertices in $G_{p_1}$ is at most $2d^{1/50}$, and so Markov's inequality shows that fewer than $n/2$ vertices are isolated, and hence that there is a non-isolated component, with probability $1-o(1)$.

Using \Cref{lem:isolates-far,lem:component-gap} we may condition on the event that every component of $G_{p_1}$ is either an isolated vertex or has order at least $d^{20}$, and that the isolated vertices are pairwise at distance at least three in $G$. We will show that after sprinkling with $p_2$, typically all the large components, i.e., the components of order at least $d^{20}$ merge.
Let $L$ be the union of the components of order at least $d^{20}$.
Consider a non-trivial partition $A\sqcup B=L$ that respects these components and, by interchanging $A$ and $B$, assume that $a\coloneqq |A|\le|B|$.
Put
\[
  A'=(A\cup N_G(A))\setminus B,
  \qquad B'=V\setminus A'.
\]
Both $A'$ and $B'$ have size at least $a$.
Indeed, they contain $A$ and $B$, respectively.
Let $D$ be the smaller of $A'$ and $B'$.
Then $a\le|D|\le n/2$, and the expansion from \eqref{eq:expansion-log} gives
\[
  e_G(A',B')=e_G(D,D^c)
  \ge \frac{e}{M}|D|\ln\frac{n}{|D|}
  \ge \alpha_Ma,
\]
where $\alpha_M$ is a constant depending solely on $M$.
For the last inequality, use the concavity of $x\mapsto x\ln(n/x)$ on $[a,n/2]$ and evaluate it at the endpoints.

For each edge of $E_G(A',B')$, choose one path of length at most three from $A$ to $B$ containing that edge.
To see this, note that $V\setminus L$ consists of isolated vertices of $G_{p_1}$, no two of which are adjacent in $G$.
An isolated vertex in $A'$ has a neighbour in $A$.
An isolated vertex in $B'$ has no neighbour in $A$; all its $G$-neighbours lie in $L$ because distinct isolated vertices are non-adjacent, and therefore it has a neighbour in $B$.
A fixed edge of $G$ belongs to at most $O(d^2)$ of the chosen paths, because the boundary edge that determines such a path is at distance at most two from it in the line graph of $G$; thus each path has at most three edges.
Greedily retaining such a path and deleting all other paths that share an edge with it therefore leaves at least
\[
  \frac{\alpha_M a}{d^2}
\]
edge-disjoint paths from $A$ to $B$ for some constant $\alpha_M>0$.
The probability that none of them is present when exposing the edges of the second round is at most
\[
  (1-p_2^3)^{\alpha_Ma/d^2}\le \exp\left\{-\frac{\alpha_Ma}{d^8}\right\}.
\]

It remains to take a union bound over the partitions of the large components in $G_{p_1}$. There are at most $N\coloneqq n/d^{20}$ large components.
For a partition whose smaller union has size $a$, that union contains at most $a/d^{20}$ components.
Consequently, the number of possible partitions is at most
\[
  \sum_{j\le a/d^{20}}\binom{N}{j}
  \le \exp\left\{O\left(\frac{a\ln n}{d^{20}}\right)\right\}.
\]
This bound follows from the standard binomial estimate when $a/d^{20}\le N/2$, and from $2^N\le\exp(2a/d^{20})$ otherwise.
Since $\ln n=O_M(d)$, summing over $d^{20}\le a\le n/2$ gives
\[
  \sum_{a=d^{20}}^{n/2}
  \exp\left\{O_M\left(\frac{a}{d^{19}}\right)-\frac{\alpha_Ma}{d^8}\right\}=o(1).
\]
Thus all non-isolated components of $G_{p_1}$ merge after exposing the edges of the second round.

Finally, an isolated vertex of $G_{p_1}$ that receives an edge in the second round cannot join another such vertex, because these vertices are not adjacent in $G$.
It therefore joins the merged component.
The vertices that receive no sprinkled edge remain isolated and by \Cref{lem:isolates-far} they are still of distance at least three.
\end{proof}

\begin{lemma}\label{lem:number-isolates}
Whp,
\[
  1\le \iso(G_{p_-})\le d
  \qquad\text{and}\qquad
  \iso(G_{p_+})=0.
\]
\end{lemma}

\begin{proof}
Let $Y_p=\iso(G_p)$.
Then $\mathbb{E}Y_p=n(1-p)^d$.
For $p=p_+$ this expectation is $1/\omega$, and Markov's inequality gives $Y_{p_+}=0$ \whp.
For $p=p_-$, we have $\mathbb{E}Y_{p_-}=\omega$, and
\[
  \mathbb{P}(Y_{p_-}>d)\le\frac{\omega}{d}=o(1).
\]

For the lower bound, expand the second moment.
Non-adjacent pairs of vertices are simultaneously isolated with probability $(\omega/n)^2$. An adjacent pair of vertices has probability at most $\alpha_M(\omega/n)^2$, where $\alpha_M$ is a suitable constant.
Since $G$ has $nd/2$ edges,
\[
  \mathbb{E}Y_{p_-}^2
  \le \omega+\omega^2+O\left(\frac{d\omega^2}{n}\right)
  =\omega^2+O(\omega).
\]
Hence $\operatorname{Var}(Y_{p_-})=O(\omega)=o(\omega^2)$, and Chebyshev's inequality gives $Y_{p_-}>0$ \whp.
\end{proof}

Combining \Cref{prop:unique-component} and \Cref{lem:number-isolates} proves part (a) and (b) of \Cref{th:window-connectivity}.



\section{A perfect matching in the random graph process}\label{s:matching}

In order to prove the last part of \Cref{th:window-connectivity} we shall prove that, with high probability, every
$G_p$ with $p\in[p_-,p_+]$ and $\delta(G_p)\ge1$ has a perfect matching.
Given the results in \Cref{s:connectivity}, every such $G_p$ is connected whp:
the isolated vertices of $G_{p_-}$ attach to the unique non-trivial component as soon as they receive an edge.

We make use of the Tutte-Berge formula. In particular, the formula (combined with the assumption that $G$ has even order) implies that the existence of a set $U \subseteq V$ with
\begin{align} \label{eq:tutte-violation}
    \odd(G_p - U) \geq |U|+2
\end{align}
is the only obstruction to a perfect matching in $G_p$.
Define
\[
 u_0\coloneqq \frac{n}{d^{M^5/p_-}}.
\]

The following two lemmas comprise our two crucial ingredients.

\begin{lemma}\label{lem:no-small-tutte}
Whp, simultaneously for every $p\in[p_-,p_+]$ with $\delta(G_p)\ge1$, the graph $G_p$ has no set $U$ with $1\le|U|\le u_0$ satisfying \eqref{eq:tutte-violation}.
\end{lemma}

\begin{lemma}\label{lem:large-ge}
Whp, simultaneously for every $p\in[p_-,p_+]$ with $\delta(G_p)\ge1$, if $\defi(G_p)\ge2$, then the set $A$ in the Gallai--Edmonds decomposition of $G_p$ has size at most $u_0$.
\end{lemma}

The proof of these two lemmas will be the task of the remaining subsections. First, we show that they indeed complete the proof of \Cref{th:window-connectivity}.

\begin{proof}
Let $p \in [p_-, p_+]$. Recall that whp if $\delta(G_p)\geq 1$, then $G_p$ is connected by part (a) and (b) of \Cref{th:window-connectivity} and we condition on this event.
Suppose that $G_p$ has no perfect matching and let $(D,A,C)$ be its Gallai--Edmonds decomposition.
Let $D_1,\ldots,D_m$ be the components of $G_p[D]$.
Put $a=|A|$ and $r=\defi(G_p)\ge2$ (since $|V(G)|$ is even).
The set $A$ is non-empty: if $A=\varnothing$, the distinct $D_i$ are components of $G_p$, contrary to connectedness and $r\ge2$.
Each $D_i$ has odd order, every component of $G_p[C]$ has even order, and there is no edge of $G_p$ between $D$ and $C$.
Consequently,
\[
 \odd(G_p-A)=m=a+r\ge a+2,
\]
so $A$ violates Tutte's condition, i.e, it satisfies \eqref{eq:tutte-violation}.
If $a\le u_0$, this contradicts \Cref{lem:no-small-tutte}.
If $a>u_0$, it contradicts \Cref{lem:large-ge}.
\end{proof}

\subsection{Tutte obstructions: Proof of \Cref{lem:no-small-tutte}}
We first state and prove several lemmas that will be of use when proving \Cref{lem:no-small-tutte}.

We start with an observation concerning the structure of minimal Tutte obstructions.
\begin{lemma} \label{lem:three_edges}
Let $U$ be a minimal set of vertices satisfying \eqref{eq:tutte-violation}. Then every vertex $v \in U$ has neighbours (in $G$) in at least three components of $G_p-U$.
\end{lemma}
\begin{proof}
Suppose towards contradiction that there is a minimal set $U$ of size $u$ satisfying \eqref{eq:tutte-violation}, with $\ell>u$ components in $G_p-U$ and containing a vertex $v$ whose neighbourhood in $G$ intersects at most two components of $G_p-U$.

Consider $U'=U \setminus \{v\}$. Clearly, $|U'|=u-1$. Observe that since $v$ is connected to at most two components in $G_p-U$, there are at least $\ell-1>u-1=|U'|$ components in $G_p-U'$. But then $U'$ is an obstruction, contradicting the minimality of $U$.
\end{proof}

Typically, subsets of vertices in $G$ expand well. For sets with small edge boundary, we use a fingerprint argument to obtain an improved bound on the number of such sets. This makes use of the Cartesian product structure of $G$.

\begin{lemma}\label{lem:small-boundary-sets}
Let $\frac{n}{d^{\ln^2 d}}\le a\le \frac{n}{2}$. Then, the number of sets $A\subseteq V$ of size $a$ with $e(A,A^c)< a\ln^2d$ is at most $\exp\left\{\frac{2a}{\ln d}\right\}$.
\end{lemma}

\begin{proof}
Let $\mathcal{F}$ be the family of $A\subseteq V$ satisfying the conditions of the lemma.

For $i\in [t]$ and any $A\subseteq V$, let $E_i(A,A^C)\subseteq E(A,A^C)$ be the set of edges in $E(A,A^C)$ corresponding to a change in the $i$-th coordinate, and let $e_{i,A}\coloneqq |E_i(A,A^C)|$. Moreover, given $I\subseteq [t]$, let $e_{I,A}=\sum_{i\in I}e_{i,A}$. We say that $A$ is \textit{bad} with respect to a set of coordinates $I$, if $e_{I,A}< a\ln^2d\cdot \frac{|I|}{t}$. Let $\mathcal{A}_{I}$ be the family of sets $A$ which are bad with respect to some $I\subseteq [t]$. Note that for every fixed $m \in \mathbb{N}$, if $A\in \mathcal{F}$, then there is some $I$ with $|I|=m$ such that $e_{I,A}< a\ln^2d\cdot \frac{|I|}{t}$. Thus,
\[|\mathcal{F}|\le \sum_{\substack{I \subseteq [t]\\|I|=m}} |\mathcal{A}_I| \leq \binom{t}{m} \max_{\substack{I \subseteq [t]\\|I|=m}} |\mathcal{A}_I|.\]
We now set $m=\log_2\left(\ln^5 d\right)$, and turn to estimate $|\mathcal{A}_I|$ for any $I \subseteq [t]$ with $|I|=m$.

For such an $I$ and $v\in V$, let $G(I,v)\coloneqq \square_{i\in [t]\setminus I}\{v_i\} \square_{i\in I}H_i\subseteq G$ be the product graph where we fix the coordinates outside of $I$ to agree with $v$. That is, $G(I, v)$ is a product graph of dimension $|I|=m$. Observe that $2^{|I|}\le|V(G(I,v))|\le M^{|I|}$, and that for every $v\neq u\in V$, $V(G(I,v))$ and $V(G(I,u))$ are either disjoint or identical. Thus, fixing $I$ with $|I|=m$, there are at most $\frac{n}{2^m}$ different subgraphs $G(I,v)$, and their union is $V$. We say that $A$ \textit{intersects non-trivially} with $G(I,v)$ if $V(G(I,v))\cap A\neq \varnothing$ and $V(G(I,v))\setminus A\neq \varnothing$.
Applying \Cref{lem:product-connectivity} to $G(I, v)$, we have that if $A$ intersects non-trivially with $G(I,v)$, then $G(I,v)$ spans at least $d(G(I, V)) \geq |I|=m$ edges of $E(A,A^C)$. Thus, if $A\in \mathcal{A}_I$, we have that $A$ intersects non-trivially at most $\frac{a\ln^2d}{t}$ such subgraphs. Indeed, otherwise, there would be no $I$ with $|I|=m$ and $e_{I,A}<a\ln^2d\cdot \frac{m}{t}$. Therefore, a set $A\in \mathcal{A}_I$ contains at least $$\frac{a-M^m\frac{a\ln^2d}{t}}{M^m} \leq \frac{a}{M^m}$$ subgraphs of the form $G(I, v)$ for $v \in V$, and at most $M^m\frac{a\ln^2d}{t}$ other vertices. We thus obtain that
\begin{align} \label{eq:bad-expansion-choices}
    |\mathcal{A}_I|&\le \binom{\frac{n}{2^m}}{\frac{a}{M^m}}\binom{n}{M^m\frac{a\ln^2d}{t}},
\end{align}
where we first choose from the at most $n/2^m$ many the subgraphs $G(I, v)$ that are contained in $A$ and then the other vertices.
We may bound \eqref{eq:bad-expansion-choices} from above using that $a \geq \frac{n}{d^{\ln^2 d}}$:
\begin{align*}
    &\left(\frac{enM^m}{a}\right)^{\frac{a}{M^m}}\left(\frac{ent}{M^m a\ln^2d}\right)^{M^m\frac{a\ln^2d}{t}}
    \le \left(ed^{\ln^2d}M^m\right)^{\frac{a}{M^m}}\left(\frac{etd^{\ln^2d}}{M^m\ln^2d}\right)^{M^m\frac{a\ln^2d}{t}}\\
    \le&\left(ed^{\ln^2d}\ln^5 d\right)^{\frac{a}{\ln^5d}}\left(\frac{ed^{\ln^2d+1}}{\ln^7d}\right)^{\frac{a\ln^7d}{t}}
    \le \exp\left\{\frac{a}{\ln d}\right\}.
\end{align*}
Altogether, by multiplying with the number of choices for $I$ we obtain that 
\[
|\mathcal{F}|\le \binom{t}{\log_M(\ln^5d)}\exp\left\{\frac{a}{\ln d}\right\}\le \exp\left\{\frac{2a}{\ln d}\right\},
\]
as required.
\end{proof}



We set
\begin{equation*}
 u_1\coloneqq \frac{n}{d^{M^3/p_-}}
\end{equation*}
and note that $u_1 \gg d\cdot u_0$.

\begin{lemma}\label{lem:no-small-separator}
Whp, there are no disjoint sets $U, B\subseteq V$ satisfying
\[
 |U|\le u_0, \qquad u_1\le|B|\le n/2,
\]
for which every edge from $B$ to $V\setminus(B\cup U)$ is not present in $G_{p_-}$.
\end{lemma}

\begin{proof}
Denote by $b \coloneqq |B|$. The number of possible $U$ is at most
\begin{align*}
    \sum_{u\le u_0}\binom nu \le\exp\left\{u_0\ln(en/u_0)\right\}=\exp\left\{u_0\cdot \frac{M^5}{p_-}\cdot \ln d \right\}\le \exp\left\{\frac{b}{\ln d}\right\},
\end{align*}
where we used the definition of $u_0$ and $u_1$ and that $b\ge u_1$.

If $e_G(B,B^c)\ge b\ln^2 d$, then at least $b\ln^2d-d|U|\ge b\ln^2 d/2$ edges from $B$ avoid $U$. Thus, the probability that there are sets $ U, B$ with $b=|B|$ as in the claim is at most
\[
 \binom nb\exp\left\{\frac{b}{\ln d}\right\}(1-p_-)^{b\ln^2 d/2}=o(n^{-3}).
\]
If $e_G(B,B^c)<b\ln^2 d$, \Cref{lem:small-boundary-sets} gives at most $\exp\{2b/\ln d\}$ choices for $B$.
Moreover the expansion property in \eqref{eq:expansion-log} gives $e_G(B,B^c)\ge \alpha_M b$ for some constant $\alpha_M$ while $du_0=o(b)$, thus at least $\frac{\alpha_M}{2} b$ edges from $B$ avoid $U$. Thus, the probability that there are sets $ U, B$ with $b=|B|$ as in the claim is at most
\[
 \exp\left\{\frac{3b}{\ln d}\right\}(1-p_-)^{\alpha_M b/2}\le\exp\{-\alpha'_M b\}.
\]
A union bound over $b \geq u_1$ concludes the proof.
\end{proof}

Before turning to the key lemma of this subsection we give bounds on the expansion of the graph that will be useful throughout the calculations.
Towards this, for $1 \leq k \leq \frac{n}{2}$ let
\begin{equation*}\label{eq:f-star}
 f^*(k)\coloneqq \max\left\{k\bigl(d-(M-1)\log_Mk\bigr),
             \frac{ek}{M}\ln\frac nk\right\}.
\end{equation*}
Extend it by $f^*(k)=f^*(n-k)$ for $n/2<k<n$.
By \Cref{th:expansion}, every non-empty proper set $A\subseteq V$ satisfies
\begin{align} \label{eq:expansion}
    e_G(A,A^c)\ge f^*(|A|).
\end{align}

\begin{lemma}\label{lem:strong-boundary}
For all sufficiently large $d$ and $3\le k\le2u_1$,
\begin{align}
 f^*(k)&\ge2d+15 (M-1) k \ln d.\label{eq:strong-3d}
\end{align}
\end{lemma}

\begin{proof}
For $k=3$, \eqref{eq:expansion} together with the first term of \eqref{eq:f-star} give
\[
 f^*(k)\ge 3d- 3(M-1) \log_M(3) \geq 2d + 45(M-1) \ln d,
\]
as claimed. For $4 \leq k \le d^3$ the same implies
\[
f^*(k)\ge kd - 3k(M-1) \log_M(d) \geq 3d + 15k (M-1) \ln d.
\]
If $d^3<k\le2u_1$, then \eqref{eq:expansion} together with the second term of \eqref{eq:f-star} give
\[
 f^*(k)\ge\frac{ek}{M}\left(\frac{M^3}{p_-}\ln d-\ln 2\right).
\]
Using that $p_-\le1-1/M$
, we have $\frac{eM^2}{p_-} \geq 16(M-1)+1$. This yields
\[
 f^*(k)\ge 16(M-1) k \ln d + k \ln d - \frac{ek}{M} \ln 2 \geq 3d + 15(M-1)k \ln d,
\]
where the last inequality uses that $k>d^3$.
\end{proof}

The following lemma is designed to show that whp there are no small sets $U$ satisfying \eqref{eq:tutte-violation}. Here, $U$ plays the role of the obstruction and $B$ contains components from $G_p-U$. Note that we treat components with an even number of vertices, except two, the same as components with an odd number of vertices, but this simplified point of view will be sufficient in this setting.


\begin{lemma}\label{lem:small-components-count}
Whp, there are no $p\in[p_-,p_+]$ and disjoint $U,B\subseteq V$ with
\[
1 \leq |U| \leq u_0, \qquad |B|<u_1
\]
such that
\begin{enumerate}[(i)]
\item no component of $B$ in $G_p-U$ has size two;
 \item the number of components of $B$ in $G_p-U$ is at least $|U|+1$;
 \item every isolated vertex of $B$ in $G_p-U$ has a $G_p$-neighbour in $U$;
 \item every vertex in $U$ has at least two neighbours in distinct components of $B$ in $G_p-U$;
 \item no edge of $G_p$ joins $B$ to $V\setminus(U\cup B)$.
\end{enumerate}
\end{lemma}

\begin{figure}
    \centering

\tikzset{every picture/.style={line width=0.75pt}} 

\begin{tikzpicture}[x=0.75pt,y=0.75pt,yscale=-1,xscale=1]

\draw   (139.33,66) .. controls (139.33,42.25) and (158.59,23) .. (182.33,23) -- (507.33,23) .. controls (531.08,23) and (550.33,42.25) .. (550.33,66) -- (550.33,195) .. controls (550.33,218.75) and (531.08,238) .. (507.33,238) -- (182.33,238) .. controls (158.59,238) and (139.33,218.75) .. (139.33,195) -- cycle ;
\draw   (164,126) .. controls (164,96.73) and (190.49,73) .. (223.17,73) .. controls (255.84,73) and (282.33,96.73) .. (282.33,126) .. controls (282.33,155.27) and (255.84,179) .. (223.17,179) .. controls (190.49,179) and (164,155.27) .. (164,126) -- cycle ;
\draw   (350,56.87) .. controls (350,43.69) and (360.69,33) .. (373.87,33) -- (445.47,33) .. controls (458.65,33) and (469.33,43.69) .. (469.33,56.87) -- (469.33,199.13) .. controls (469.33,212.31) and (458.65,223) .. (445.47,223) -- (373.87,223) .. controls (360.69,223) and (350,212.31) .. (350,199.13) -- cycle ;
\draw  [fill={rgb, 255:red, 0; green, 0; blue, 0 }  ,fill opacity=1 ] (386,60) .. controls (386,57.79) and (387.79,56) .. (390,56) .. controls (392.21,56) and (394,57.79) .. (394,60) .. controls (394,62.21) and (392.21,64) .. (390,64) .. controls (387.79,64) and (386,62.21) .. (386,60) -- cycle ;
\draw  [fill={rgb, 255:red, 0; green, 0; blue, 0 }  ,fill opacity=1 ] (432,92) .. controls (432,89.79) and (433.79,88) .. (436,88) .. controls (438.21,88) and (440,89.79) .. (440,92) .. controls (440,94.21) and (438.21,96) .. (436,96) .. controls (433.79,96) and (432,94.21) .. (432,92) -- cycle ;
\draw  [fill={rgb, 255:red, 0; green, 0; blue, 0 }  ,fill opacity=1 ] (402,102) .. controls (402,99.79) and (403.79,98) .. (406,98) .. controls (408.21,98) and (410,99.79) .. (410,102) .. controls (410,104.21) and (408.21,106) .. (406,106) .. controls (403.79,106) and (402,104.21) .. (402,102) -- cycle ;
\draw  [fill={rgb, 255:red, 0; green, 0; blue, 0 }  ,fill opacity=1 ] (416,73) .. controls (416,70.79) and (417.79,69) .. (420,69) .. controls (422.21,69) and (424,70.79) .. (424,73) .. controls (424,75.21) and (422.21,77) .. (420,77) .. controls (417.79,77) and (416,75.21) .. (416,73) -- cycle ;
\draw   (361,133) .. controls (361,126.37) and (370.48,121) .. (382.17,121) .. controls (393.86,121) and (403.33,126.37) .. (403.33,133) .. controls (403.33,139.63) and (393.86,145) .. (382.17,145) .. controls (370.48,145) and (361,139.63) .. (361,133) -- cycle ;
\draw   (372,163) .. controls (372,156.37) and (381.48,151) .. (393.17,151) .. controls (404.86,151) and (414.33,156.37) .. (414.33,163) .. controls (414.33,169.63) and (404.86,175) .. (393.17,175) .. controls (381.48,175) and (372,169.63) .. (372,163) -- cycle ;
\draw   (411,184) .. controls (411,177.37) and (420.48,172) .. (432.17,172) .. controls (443.86,172) and (453.33,177.37) .. (453.33,184) .. controls (453.33,190.63) and (443.86,196) .. (432.17,196) .. controls (420.48,196) and (411,190.63) .. (411,184) -- cycle ;
\draw   (409.67,132.5) .. controls (409.67,125.87) and (419.14,120.5) .. (430.83,120.5) .. controls (442.52,120.5) and (452,125.87) .. (452,132.5) .. controls (452,139.13) and (442.52,144.5) .. (430.83,144.5) .. controls (419.14,144.5) and (409.67,139.13) .. (409.67,132.5) -- cycle ;
\draw [color={rgb, 255:red, 74; green, 144; blue, 226 }  ,draw opacity=1 ]   (207.33,84) .. controls (251.33,57) and (356.33,51) .. (390,60) ;
\draw [color={rgb, 255:red, 74; green, 144; blue, 226 }  ,draw opacity=1 ]   (262.33,111) .. controls (313.33,92) and (375.33,93) .. (406,102) ;
\draw [color={rgb, 255:red, 74; green, 144; blue, 226 }  ,draw opacity=1 ]   (228.33,91) .. controls (274.33,65) and (385.33,68) .. (420,73) ;
\draw [color={rgb, 255:red, 74; green, 144; blue, 226 }  ,draw opacity=1 ]   (244.33,101) .. controls (294.33,80) and (399.33,83) .. (436,92) ;
\draw [color={rgb, 255:red, 208; green, 2; blue, 27 }  ,draw opacity=1 ][line width=3]  [dash pattern={on 7.88pt off 4.5pt}]  (443.33,130) -- (512.33,125) ;
\draw [color={rgb, 255:red, 208; green, 2; blue, 27 }  ,draw opacity=1 ][line width=3]  [dash pattern={on 7.88pt off 4.5pt}]  (436,92) -- (505,87) ;
\draw [color={rgb, 255:red, 208; green, 2; blue, 27 }  ,draw opacity=1 ][line width=3]  [dash pattern={on 7.88pt off 4.5pt}]  (442.33,187) -- (506.33,213) ;
\draw  [fill={rgb, 255:red, 0; green, 0; blue, 0 }  ,fill opacity=1 ] (242,137) .. controls (242,134.79) and (243.79,133) .. (246,133) .. controls (248.21,133) and (250,134.79) .. (250,137) .. controls (250,139.21) and (248.21,141) .. (246,141) .. controls (243.79,141) and (242,139.21) .. (242,137) -- cycle ;
\draw    (246,137) -- (372.33,134) ;
\draw    (246,137) -- (406,102) ;
\draw    (244,156) -- (375.33,139) ;
\draw  [fill={rgb, 255:red, 0; green, 0; blue, 0 }  ,fill opacity=1 ] (240,156) .. controls (240,153.79) and (241.79,152) .. (244,152) .. controls (246.21,152) and (248,153.79) .. (248,156) .. controls (248,158.21) and (246.21,160) .. (244,160) .. controls (241.79,160) and (240,158.21) .. (240,156) -- cycle ;
\draw    (244,156) -- (384.33,165) ;

\draw (147,77) node [anchor=north west][inner sep=0.75pt]   [align=left] {$\displaystyle U$};
\draw (482,39) node [anchor=north west][inner sep=0.75pt]   [align=left] {$\displaystyle B$};
\draw (403,40) node [anchor=north west][inner sep=0.75pt]   [align=left] {$\displaystyle Q_{1}$};
\draw (400,195) node [anchor=north west][inner sep=0.75pt]   [align=left] {$\displaystyle Q_{\ell +r}{}$};

\end{tikzpicture}

    \caption{Sets $U, B$ as in \Cref{lem:small-components-count}. All components of $B$ in $G_p-U$ are either isolated vertices or of size at least three. Isolated vertices of $B$ in $G_p-U$ have a $G_p$-neighbour in $U$ (blue edges) and vertices in $U$ have  $G_p$-neighbours in at least two distinct components of $B$ in $G_p-U$ (black edges). The dashed red edges are not present in $G_p$.}
    \label{fig:small-components-count}
\end{figure}

\begin{proof}
Denote by $u\coloneqq |U|$ and by $b\coloneqq |B|$. Let $Q_1,...,Q_{\ell}$ be the components of $B$ containing a single isolated vertex in $G_p-U$, and let $Q_{\ell+1},..., Q_{\ell+r}$ be the components of $B$ size at least three in $G_p-U$ (see \Cref{fig:small-components-count}). Then $\ell+r$ is the number of components of $B$ in $G_p-U$. Further, we set $s=\sum_{i=\ell+1}^{\ell+r} |Q_i|$ to be the number of vertices contained in components of size at least three.
On $U\cup B$, define a subgraph $J$ of $G$ by
\[
 E(J)=E_G(B)\cup E_G(U,B).
\]
By (iii) and (iv) every component $K_1,\ldots,K_a$ of $J$ in $G$ has order at least three, and hence $a\le \frac{b}{3}$. 
By assumption all edges between $Q_i$ and $Q_j$ for $i, j \leq \ell+r$ with $i \neq j$ are not present in $G_{p_-}$. Furthermore, by (iv), any edge from $B$ to $V\setminus (U \cup B)$ is not present in $G_{p_-}$. Denote by $h$ the total number of these edges. Then,
\begin{equation*}\label{eq:small-components-key}
 2h \ge d\ell+\sum_{i=\ell+1}^r e_G(Q_i,Q_i^c)+\sum_{j=1}^c e_G(K_j,K_j^c)-du.
\end{equation*}
Indeed, every isolated vertex from a component $Q_i$ for $i \leq \ell$ has degree $d$ in $G$ and every other component $Q_i$ for $\ell+1 \leq i \leq \ell+r$ contributes $e_G(Q_i, Q_i^c)$ edges. With this, we have counted each edge between $Q_i$ and $Q_j$ for $i,j\leq \ell+r$ with $i\neq j$ exactly twice, and the edges between $B$ and $V\setminus (U\cup B)$ at most once. 
Adding the edges from $e_G(K_i, K_i^c)$ counts each edge at most twice. In the first two terms, we may have also counted edges between $U$ and $B$ (each such edge at most once), and in the second term we may have counted edges spanned by $U$ (at most twice); all of these are accounted for when we subtract  $du$.  

Now, we use \Cref{lem:strong-boundary} (noting that all the components under consideration are of size at most $2u_1$) to give a lower bound on the number of absent edges. Using the expansion from \eqref{eq:strong-3d}, we obtain
\[
2h \geq d(\ell-u) + 2dr + 15 s (M-1) \ln d+2da+ 15 (b+u)(M-1) \ln d.
\]
Now $\ell+r \geq u+1$, so $\ell-u+2r \geq 1$ and consequently
\[
h \geq da + \frac{d}{2}+\frac{15(M-1)}{2} (b+u+s) \ln d.
\]

It remains to bound the number of choices for $U, B$ satisfying the conditions. Fix the number $a$ of components of $J$ and $b+u$. First, we choose a root for each component and then a rooted spanning tree marking each of its vertices as belonging to $U$ or $B$. These choices determine $U, B$ and the component structure of $B$ in $G_p-U$ and using \Cref{lem:rooted-trees} there are at most
\[
n^a (ed)^{b+u} 2^{b+u} \leq n^a d^{2(b+u)}
\]
possibilities.
Combining with a union bound over the number of choices for $b+u$ and $a$, we conclude that the probability to have sets $B, U$ as in the claim is at most
\begin{align} \label{eq:probability-small-Tutte}
    \sum_{b+u=3}^{2u_1} \sum_{a \leq (b+u)/3} n^a d^{2(b+u)} (1-p_-)^{da + \frac{d}{2}+\frac{15(M-1)}{2} (b+u) \ln d}.
\end{align}
Using that $(1-p_-)^d= \frac{\omega}{n}$ and $\ln n \geq \frac{\ln 2}{M-1} d$, the probability in \eqref{eq:probability-small-Tutte} is bounded by
\[
\left(\frac{\omega}{n} \right)^{1/2} \sum_{j=3}^{2 u_1} j \left(d^2 \omega^{1/3} \left( \frac{\omega}{n}\right)^{\frac{15(M-1)}{2} \frac{\ln d}{d}} \right)^j = \left(\frac{\omega}{n} \right)^{1/2} \sum_{j\geq 3} j \exp\left( - \Omega(j \ln d)\right) = o(1),
\]
since $\frac{15 \ln 2}{2} \geq 5$.

\end{proof}

We are now ready to derive Lemma \ref{lem:no-small-tutte}.
\begin{proof}[Proof of \Cref{lem:no-small-tutte}]
Recall that whp every $G_p$ under consideration is connected and that $G$ has even order. In particular, the empty set does not satisfy \eqref{eq:tutte-violation}.
Choose a set $U$ satisfying \eqref{eq:tutte-violation} of minimum cardinality, put $u=|U|$, and let $K_1,\ldots,K_m$ be the odd components of $G_p-U$.
By \Cref{lem:three_edges} every $v\in U$ has neighbours in at least three of these components. 
Every $K_i$ which is an isolated vertex in $G_p-U$ has a $G_p$-neighbour in $U$, because $\delta(G_p)\ge1$.

Suppose that $u \leq u_0$.  
Whp, at most one $K_i$ has order larger than $u_1$ by \Cref{lem:no-small-separator}. Let $B$ be the union of the remaining odd components. Note that there are at least $u+1$ of them and every $v\in U$ has neighbours in at least two components of $B$.
If $|B| < u_1$, then $U$ and $B$ satisfy the assumptions of \Cref{lem:small-components-count}, a contradiction.
If $|B|\geq u_1$, then a greedy union of its components has size between $u_1$ and $2u_1<n/2$. In this case $U$ together with this union of components satisfy the assumptions of \Cref{lem:no-small-separator}, a contradiction.
\end{proof}

\subsection{Edmonds-Gallai decomposition: Proof of \Cref{lem:large-ge}}

Fix $p\in[p_-,p_+]$ with $\delta(G_p)\ge1$, and suppose that $G_p$ has no perfect matching. Recall, that whp $G_p$ is connected. 
We use the notation set out in \Cref{s:Gallai} and abbreviate
\[
 D=D_{G_p},\qquad A=A_{G_p},\qquad C=C_{G_p}.
\]
Let $D_1,\ldots,D_m$ be the components of $G_p[D]$.
Set
\[
 a=|A|,\qquad r=m-a=\defi(G_p)\ge2.
\]
Here $r\ge2$ because $G_p$ has even order and no perfect matching.
We start with some structural observations about the Edmonds-Gallai decomposition.

Let $Q$ be the auxiliary bipartite graph with parts $A$ and $\{D_1,\ldots,D_m\}$, where $v\in A$ is adjacent to $D_i$ whenever $v$ has a $G_p$-neighbour in $D_i$.

\begin{lemma}\label{lem:strict-hall}
For every non-empty $S\subseteq A$,
\[
 |N_Q(S)|\ge|S|+1.
\]
Consequently every vertex of $A$ has neighbours in at least two components of $G_p[D]$, and $Q$ has at most $r$ connected components.
\end{lemma}

\begin{proof}
For every $i$ and every $v\in D_i$, the definition of $D$ gives a maximum matching $M_v$ leaving $v$ unmatched; in particular, in $M_v$ there is no edge from $A$ to $D_i$. 
The edges in $M_v$ match $A$ to distinct members of $\{D_j:j\ne i\}$.
Thus $Q-\{D_i\}$ has a matching saturating $A$.
In particular, $Q$ has a matching saturating $A$, so Hall's theorem gives $|N_Q(S)|\ge|S|$.
If $|N_Q(S)|=|S|$, deleting any $D_i\in N_Q(S)$ would violate Hall's condition, a contradiction.

In each component of $Q$, the number of components of $G_p[D]$ is therefore at least one more than the number of vertices in $A$.
Summed over all components of $Q$, this difference is $m-a=r$, so $Q$ has at most $r$ components.
\end{proof}

Let
\[
 h\coloneqq \sum_{i<j}e_G(D_i,D_j)+e_G(D,C).
\]
and note that all edges counted in $h$ are not present in $G_p$ (and each edge is counted once).
We start by giving a lower bound on $h$.

\begin{lemma}\label{lem:ge-boundary-identity}
\begin{equation}\label{eq:ge-exact}
 2h=dr+2e_G(A)+e_G(A\cup D,C)  +\sum_{i:|D_i|\ge3}\bigl(e_G(D_i,D_i^c)-d\bigr).
\end{equation}
\end{lemma}

\begin{proof}
Degree summation in $D$ gives
\[
    d |D|= 2 \sum_{i} e_G(D_i)+ 2 \sum_{i<j} e_G(D_i, D_j) + e_G(D, A) + e_G(D, C).
\]
Rearranging for $h=\sum_{i<j}e_G(D_i,D_j)+e_G(D,C)$, we obtain
\begin{align} \label{eq:summation-D}
    2h=d|D|-2\sum_ie_G(D_i)-e_G(D,A)+e_G(D,C).
\end{align}
Similarly, degree summation in $A$ gives
\begin{align} \label{eq:summation-A}
    d |A| = 2 e_G(A) + e_G(A, D) + e_G(A, C).
\end{align}
Subtracting the LHS and adding the RHS of \eqref{eq:summation-A} to \eqref{eq:summation-D}, we obtain
\[
2h= d( |D|-|A|) - 2 \sum_{i} e_G(D_i) + 2 e_G(A) + e_G(A, C) + e_G(D, C).
\]
Now
\[
 |D|-|A|=r+\sum_{|D_i|\ge3}(|D_i|-1),  \qquad 2e_G(D_i)=d|D_i|-e_G(D_i,D_i^c),
\]
and since each $|D_i|$ is odd both sums over $i$ may be combined, proving the claim.
\end{proof}

We first deal with the case where $G$ is bipartite and $G_p[D]$ contains only isolated vertices.
\begin{lemma}\label{lem:bipartite-singletons}
Suppose that $G$ is bipartite.
Whp, for every $p\in[p_-,p_+]$ with $\delta(G_p)\ge1$, there is no Gallai--Edmonds decomposition of $G_p$ for which $\defi(G_p)\ge2$ and every component of $G_p[D]$ is an isolated vertex.
\end{lemma}

\begin{figure}
    \centering

 
\tikzset{
pattern size/.store in=\mcSize, 
pattern size = 5pt,
pattern thickness/.store in=\mcThickness, 
pattern thickness = 0.3pt,
pattern radius/.store in=\mcRadius, 
pattern radius = 1pt}
\makeatletter
\pgfutil@ifundefined{pgf@pattern@name@_8t1jihue9}{
\pgfdeclarepatternformonly[\mcThickness,\mcSize]{_8t1jihue9}
{\pgfqpoint{0pt}{0pt}}
{\pgfpoint{\mcSize+\mcThickness}{\mcSize+\mcThickness}}
{\pgfpoint{\mcSize}{\mcSize}}
{
\pgfsetcolor{\tikz@pattern@color}
\pgfsetlinewidth{\mcThickness}
\pgfpathmoveto{\pgfqpoint{0pt}{0pt}}
\pgfpathlineto{\pgfpoint{\mcSize+\mcThickness}{\mcSize+\mcThickness}}
\pgfusepath{stroke}
}}
\makeatother

 
\tikzset{
pattern size/.store in=\mcSize, 
pattern size = 5pt,
pattern thickness/.store in=\mcThickness, 
pattern thickness = 0.3pt,
pattern radius/.store in=\mcRadius, 
pattern radius = 1pt}
\makeatletter
\pgfutil@ifundefined{pgf@pattern@name@_buu9o8d8i}{
\pgfdeclarepatternformonly[\mcThickness,\mcSize]{_buu9o8d8i}
{\pgfqpoint{0pt}{0pt}}
{\pgfpoint{\mcSize+\mcThickness}{\mcSize+\mcThickness}}
{\pgfpoint{\mcSize}{\mcSize}}
{
\pgfsetcolor{\tikz@pattern@color}
\pgfsetlinewidth{\mcThickness}
\pgfpathmoveto{\pgfqpoint{0pt}{0pt}}
\pgfpathlineto{\pgfpoint{\mcSize+\mcThickness}{\mcSize+\mcThickness}}
\pgfusepath{stroke}
}}
\makeatother

 
\tikzset{
pattern size/.store in=\mcSize, 
pattern size = 5pt,
pattern thickness/.store in=\mcThickness, 
pattern thickness = 0.3pt,
pattern radius/.store in=\mcRadius, 
pattern radius = 1pt}
\makeatletter
\pgfutil@ifundefined{pgf@pattern@name@_zi8dgom70}{
\pgfdeclarepatternformonly[\mcThickness,\mcSize]{_zi8dgom70}
{\pgfqpoint{0pt}{0pt}}
{\pgfpoint{\mcSize+\mcThickness}{\mcSize+\mcThickness}}
{\pgfpoint{\mcSize}{\mcSize}}
{
\pgfsetcolor{\tikz@pattern@color}
\pgfsetlinewidth{\mcThickness}
\pgfpathmoveto{\pgfqpoint{0pt}{0pt}}
\pgfpathlineto{\pgfpoint{\mcSize+\mcThickness}{\mcSize+\mcThickness}}
\pgfusepath{stroke}
}}
\makeatother
\tikzset{every picture/.style={line width=0.75pt}} 

\begin{tikzpicture}[x=0.75pt,y=0.75pt,yscale=-1,xscale=1]

\draw    (204.33,185.67) -- (300.33,186.67) ;
\draw    (357.33,185.67) -- (453.33,186.67) ;
\draw    (203.33,104.67) -- (301.33,104.67) ;
\draw    (356.33,104.67) -- (454.33,104.67) ;
\draw  [draw opacity=0][fill={rgb, 255:red, 65; green, 117; blue, 5 }  ,fill opacity=0.55 ] (204.33,104.17) -- (302.33,104.17) -- (302.33,185.67) -- (204.33,185.67) -- cycle ;
\draw   (204.33,22.67) -- (301.33,22.67) -- (301.33,275.67) -- (204.33,275.67) -- cycle ;
\draw   (356.33,23) -- (453.33,23) -- (453.33,276) -- (356.33,276) -- cycle ;
\draw  [draw opacity=0][fill={rgb, 255:red, 65; green, 117; blue, 5 }  ,fill opacity=0.55 ] (355.33,23) -- (453.33,23) -- (453.33,104.67) -- (355.33,104.67) -- cycle ;
\draw  [draw opacity=0][pattern=_8t1jihue9,pattern size=6pt,pattern thickness=0.75pt,pattern radius=0pt, pattern color={rgb, 255:red, 0; green, 0; blue, 0}] (204.33,186.17) -- (302.33,186.17) -- (302.33,275.67) -- (204.33,275.67) -- cycle ;
\draw  [draw opacity=0][pattern=_buu9o8d8i,pattern size=6pt,pattern thickness=0.75pt,pattern radius=0pt, pattern color={rgb, 255:red, 0; green, 0; blue, 0}] (204.33,22.67) -- (302.33,22.67) -- (302.33,104.17) -- (204.33,104.17) -- cycle ;
\draw  [draw opacity=0][pattern=_zi8dgom70,pattern size=6pt,pattern thickness=0.75pt,pattern radius=0pt, pattern color={rgb, 255:red, 0; green, 0; blue, 0}] (356.33,23.17) -- (454.33,23.17) -- (454.33,104.67) -- (356.33,104.67) -- cycle ;
\draw [color={rgb, 255:red, 208; green, 2; blue, 27 }  ,draw opacity=1 ][line width=1.5]  [dash pattern={on 5.63pt off 4.5pt}]  (292.33,45.67) -- (368.33,46.67) ;
\draw [color={rgb, 255:red, 208; green, 2; blue, 27 }  ,draw opacity=1 ][line width=1.5]  [dash pattern={on 5.63pt off 4.5pt}]  (289.33,68.67) -- (365.33,69.67) ;
\draw [color={rgb, 255:red, 208; green, 2; blue, 27 }  ,draw opacity=1 ][line width=1.5]  [dash pattern={on 5.63pt off 4.5pt}]  (287.33,227.67) -- (367.33,83.67) ;

\draw (176,10) node [anchor=north west][inner sep=0.75pt]   [align=left] {$\displaystyle L$};
\draw (482,12) node [anchor=north west][inner sep=0.75pt]   [align=left] {$\displaystyle R$};
\draw (459,62) node [anchor=north west][inner sep=0.75pt]   [align=left] {$\displaystyle D\ \cap \ R$};
\draw (459,142) node [anchor=north west][inner sep=0.75pt]   [align=left] {$\displaystyle A\ \cap \ R$};
\draw (460,227) node [anchor=north west][inner sep=0.75pt]   [align=left] {$\displaystyle C \ \cap \ R$};
\draw (147,58) node [anchor=north west][inner sep=0.75pt]   [align=left] {$\displaystyle D\ \cap \ L$};
\draw (147,135) node [anchor=north west][inner sep=0.75pt]   [align=left] {$\displaystyle A\ \cap \ L$};
\draw (140,218) node [anchor=north west][inner sep=0.75pt]   [align=left] {$\displaystyle C \ \cap \ L$};

\end{tikzpicture}

    \caption{The setup for \Cref{lem:bipartite-singletons}. The set $S$ is depicted in green while the set $D^*$ is depicted in stripes. The dashed red edges are not present in $G_p$.}
    \label{fig:bipartite-singletons}
\end{figure}

\begin{proof}
Let $V=L\sqcup R$ be the bipartition of $G$ and note that since $G$ is regular, $|L|=|R|$. Since $G_p[C]$ has a perfect matching,
\[
 |C\cap L|=|C\cap R|.
\]
Let
\[
 S=(A\cap L)\cup(D\cap R), \qquad D^*=D\cup(C\cap L),
\]
and suppose w.l.o.g. that $|(A\cap L)\cup(D\cap R)| \leq |(A\cap R)\cup(D\cap L)|$ (otherwise set $S=(A \cap R) \cup (D \cap L)$ and $D^*=D \cup (C \cap R)$).
Since $L$ is an independent set in $G$, naturally so is $C\cap L$. Further, since there are no edges between distinct components of $G_p[D]$ and we assume all components of $G_p[D]$ are singletons, there are no edges in $G_p[D]$. Finally, there are no edges in $G_p$ between $D$ and $C$ by construction. Thus, every edge of $G[D^*]$ is not present in $G_p$. 
Moreover,
\begin{align} \label{eq:def-D*}
     |D^*|-|V\setminus D^*|=r.
\end{align}
On the other hand,
\[
 e_G(D^*)+e_G(V\setminus D^*)=e_G(S,S^c),
\]
see \Cref{fig:bipartite-singletons}.
Regularity and \eqref{eq:def-D*} also give
\[
 2\bigl(e_G(D^*)-e_G(V\setminus D^*)\bigr)=dr.
\]
Combining these two equations gives the identity
\begin{equation}\label{eq:bipartite-flip}
 e_G(D^*)=\frac{dr}{4}+\frac12e_G(S,S^c).
\end{equation}
Consider the auxiliary graph $Q$ and note that by \Cref{lem:strict-hall} every vertex in $A \cap L$ has neighbours in at least two components of $D\cap R$.
Furthermore, since $G_p$ is whp connected every component of $D\cap R$ is whp connected to at least one vertex from $A \cap L$. We assume this holds deterministically, noting that this intersection changes the error probability only by a multiplicative factor $1+o(1)$. Thus, every component of $S$ in $G_p$ has order at least three. Now, we turn to bounding the number of components $\ell$ of $S$ in $G_p$. Note, that from the fact that every vertex in $A \cap L$ connects at least two components of $D\cap R$, we obtain
\[
\ell \leq |D \cap R|-|A \cap L|.
\]
Now, since $r=|D \cap L| + |D \cap R| - |A \cap L| - |A \cap R|$ and our assumption that $|(A\cap L)\cup(D\cap R)| \leq |(A\cap R)\cup(D\cap L)|$ we conclude that $\ell \leq r/2$, i.e., that there are at most $r/2$ connected components in $S$. By \eqref{eq:bipartite-flip}, at least
\[
\frac{dr}{4} + \frac{1}{2} e_G(S, S^c)
\]
edges are not present in $G_p$. We distinguish a few cases according to $S$ and in each case we bound the probability of the event in the claim by $(1-p_-)^{\frac{dr}{4} + \frac{1}{2} e_G(S, S^c)}$ and union bound over the number of choices for $S$.

If $|S|\eqqcolon s\le2u_1$, then the expansion from \eqref{eq:strong-3d} and the fact that each component has order at least three give
\[
 e_G(S,S^c)\ge 2d\ell+15(M-1) s \ln d.
\]
By \Cref{lem:rooted-trees}, there are most $n^\ell(ed)^s$ choices for $S$.
Thus the probability to have such a decomposition is at most
\begin{align} \label{eq:prob-S-Gallai}
    \sum_{r\ge2}\sum_{\ell\le r/2}\sum_{3\ell\le s\le2u_1} n^\ell(ed)^s (1-p_-)^{\frac{dr}{4}+dc+7(M-1) s \ln d}.
\end{align}
Using that $(1-p_-)^d=\frac{\omega}{n}$ and $(1-p_-)^{7(M-1) s \ln d} \leq d^{-2s}$ (since $\ln n \geq \frac{\ln 2}{M-1} d$ and $7 \ln 2 \geq 2$), the probability in \eqref{eq:prob-S-Gallai} is bounded by
\[
\sum_{r\ge2}\sum_{c\le r/2}\sum_{3c\le s\le2u_1} \omega^c e^s d^{-s} \left(\frac{\omega}{n}\right)^{\frac{r}{4}} \leq \sum_{r \geq 2} r \left( \frac{\omega^3}{n}\right)^{\frac{r}{4}} = o(1),
\]
where we used that the sum over $s$ is geometric and $c \leq r/2$.

If $s>2u_1$ and if $e_G(S,S^c)\ge s\ln^2 d$, then
\[
 \sum_s \binom{n}{s} (1-p_-)^{s\ln^2 d/2}=o(1),
\]
since $s>2u_1$ gives $\ln(en/s)=O_M(\ln d)$.

If $s>2u_1$ and if $e_G(S,S^c)< s\ln^2 d$, then by \Cref{lem:small-boundary-sets} there are at most $\exp\{2s/\ln d\}$ choices for $S$.
On the other hand, \eqref{eq:expansion-log} gives
\[
 e_G(S,S^c)\ge\frac{e\ln 2}{M}s,
\]
and hence
\[
 \sum_{s>2u_1}\exp\{2s/\ln d\}(1-p_-)^{e_G(S,S^c)/2}=o(1).
\]
\end{proof}

Next, we treat products with a non-bipartite base graph.
We shall use the following counting lemma for bipartitions of vertex sets with few edges inside the bipartition classes. The two bipartition classes are labelled and either one may be empty.

\begin{lemma}
\label{lem:two-colour-count}
Fix constants $K, K'>0$.
Let $R\subseteq V(G)$ have size $|R|\ge d$ and
\[
 e_G(R,R^c)\le K |R|\ln d.
\]
The number of bipartitions $R=R_1\sqcup R_2$ satisfying
\[
 e_G(R_1)+e_G(R_2)\le K' |R|
\]
is at most
\[
 \exp\left\{O\left(\frac{|R|\ln^2 d}{d}\right)\right\}.
\]
Consequently, for every fixed $\theta>0$,
\begin{equation}\label{eq:colouring-weight}
 \sum_{R=R_1\sqcup R_2} (1-p_-)^{\theta(e_G(R_1)+e_G(R_2))}
 \le\exp\{-o(|R|)\}.
\end{equation}
\end{lemma}

\begin{proof}
Put $H=G[R]$ and $m=e_G(R_1)+e_G(R_2)$, and assume w.l.o.g. that $|R_2|\ge|R_1|$.
Note that
\begin{align} \label{eq:counting-lem}
   e_G(R, R^c)=d|R|- \sum_{v \in R} d_H(v), \qquad 2m= \sum_{v \in R_1} d_{H[R_1]} (v) + \sum_{v \in R_2} d_{H[R_2]}(v). 
\end{align}
Let $X$ be the set of vertices $v \in R$ with either few neighbours in $H$ or many neighbours in their own bipartition class, i.e.,
\[
 X\coloneqq \left\{v \in R: d_H(v)<\frac{3}{4}d \qquad\text{or}\qquad  d_{H[R_i]}(v)>d/4\quad( \text{ if } v\in R_i)\right\}.
\]
By definition every $v \in R\setminus V$ has at least $\frac{d}{2}$ neighbours in the opposite bipartition class.
From \eqref{eq:counting-lem} we conclude that
\begin{align} \label{eq:size-X}
    |X| \leq \frac{4 e_G(R, R^c)}{d} + \frac{8m}{d}.
\end{align}

Form a random set $T$ by choosing each vertex of $R_2$ independently with probability $\rho=12\ln d/d$. The expected number of vertices of $R_1\setminus X$ with no neighbour in $T$ is at most
\[
 s(1-\rho)^{d/2}\le \frac{|R|}{d^6}.
\]
Moreover,
\[
 \mathbb{E}|T|= \rho |R_2| \le\rho |R|,
 \qquad
 \mathbb{E}|R_2\cap N_H(T)|\le\rho 2 e_G(R_2)\le2\rho m.
\]
Applying Markov's inequality to these three quantities gives a set $T\subseteq R_2$ with
\[
 |T|\le3\rho |R|,
 \qquad
 |R_1\setminus N_H(T)|\le |X| + \frac{|R|}{d^2},
 \qquad
 |R_2\cap N_H(T)|\le6\rho m.
\]
Combining \eqref{eq:size-X} with our assumptions on $e_G(R, R^c)$ and $m$, we obtain that $R_1\mathbin{\triangle}N_H(T)$ has size $O(|R|\ln d/d)$.
The pair $(T, R_1\mathbin{\triangle}N_H(T))$ determines $R_1$, and the number of possible pairs is at most
\[
 \left(\sum_{j \le O(|R|\ln d/d)}\binom{|R|}{j}\right)^2
 =\exp\left\{O\left(\frac{|R|\ln^2 d}{d}\right)\right\},
\]
which proves the first assertion.

For \eqref{eq:colouring-weight} apply the first assertion for the bipartitions with $m\le \frac{2(M-1)}{\Theta} |R|$ to obtain that these contribute at most $\exp\{o(|R |)\}$. For the bipartitions with $m>\frac{2(M-1)}{ \Theta} |R|$ contribute at most
\[
 2^{|R|} (1-p_-)^{2 (M-1) |R|}=\exp\left(- \Omega(|R|) \right),
\]
using again that $\ln n \geq \frac{\ln 2}{M-1} d$.
\end{proof}
We can now show that if $G$ is non-bipartite, then the set $A$ in the Gallai-Edmonds decomposition of $G_p$ cannot be too big.
\begin{lemma}\label{lem:nonbipartite-ge}
Suppose that at least one base graph of $G$ is non-bipartite.
Whp for every $p\in[p_-,p_+]$, if $\defi(G_p)\ge2$, then the set $A$ in the Gallai--Edmonds decomposition of $G_p$ has size at most $u_0$.
\end{lemma}
\begin{proof}
Fix $p\in[p_-,p_+]$, let $(D,A,C)$ be the Gallai--Edmonds decomposition of $G_p$, and suppose that $\defi(G_p)\ge2$ and $|A|>u_0$.
Let $D_0$ be the union of the singleton components of $G_p[D]$, let $D^+$ be the family of non-singleton components of $G_p[D]$, and put
\[
 D_1=\bigcup_{D_i \in D^+} D_i.
\]
Write $a=|A|$, $r=|D^+|+|D_0|-a\ge2$, and recall
\[
 h=\sum_{i<j}e_G(D_i,D_j)+e_G(D,C).
\]
All the edges counted in $h$ are not present in $G_{p_-}$.
Hence, we shall sum $(1-p_-)^h$ over the choices of $A$, $C$, the singleton components of $D$, and the non-singleton components of $D$.

Every $D_i\in D^+$ is a non-singleton factor-critical component, so $|D_i|\ge3$.
By \eqref{eq:ge-exact},
\begin{equation}\label{eq:mixed-exact}
 h=\frac{dr}{2}+e_G(A)+\frac{e_G(C, C^c)}{2}
   +\frac{1}{2}\sum_{D_i\in D^+}\bigl(e_G(D_i, D_i^c)-d\bigr).
\end{equation}
Also 
\begin{align} \label{eq:h-lower-bound}
    h\ge e_G(D_0),
\end{align}
since every edge of $G[D_0]$ joins two distinct components of $G_p[D]$.
Taking three quarters of \eqref{eq:mixed-exact} and one quarter of \eqref{eq:h-lower-bound} gives
\begin{align} \label{eq:combined-LB-h}
     h\ge\frac{3dr}{8}+\frac34e_G(A)+\frac14e_G(D_0)  +\frac38 e_G(C, C^c)
 +\frac38\sum_{D_i\in D^+}\bigl(e_G(D_i, D_i^c)-d\bigr).
\end{align}

Put $k_i=\min\{|D_i|,|D_i^c|\}$ so that $3\le k_i\le n/2$.
If $k_i\le2u_1$, then \eqref{eq:strong-3d} gives $e_G(D_i, D_i^c)\ge2d$.
If $k_i>2u_1$, then $\ln(n/k_i)\ge\ln 2$, and \eqref{eq:expansion-log} gives
\[
 e_G(D_i, D_i^c)\ge\frac{e\ln 2}{M}k_i>2d
\]
since $u_1\gg d$. Thus,
\begin{align} \label{eq:apply-D_i}
    \frac18\bigl(e_G(D_i, D_i^c)-d\bigr)\ge\frac1{16}e_G(D_i, D_i^c).
\end{align}
Furthermore, 
 \begin{align} \label{eq:replace-Y}
     e_G(C, C^c)+\sum_{D_i\in D^+}e_G(D_i, D_i^c) \geq e_G((C \cup D_1), (C \cup D_1)^c).
 \end{align}
Plugging \eqref{eq:apply-D_i} and \eqref{eq:replace-Y} into \eqref{eq:combined-LB-h} we obtain the following lower bound for $h$:
\begin{equation}\label{eq:mixed-master}
 h\ge\frac{3dr}{8}
 +\frac{e_G(A)+e_G(D_0)}4
 +\frac{e_G((C \cup D_1), (C, \cup D_1)^c)}{16}
 +\frac14\sum_{D_i \in D^+}\bigl(e_G(D_i, D_i^c)-d\bigr).
\end{equation}

We first sum over the possible families $D^+$ for a fixed set $(C \cup D_1)$.
By \Cref{lem:rooted-trees} here are at most $|C \cup D_1|(ed)^{k-1}$ connected sets of cardinality $k$ with a root in $C \cup D_1$.

For $3\le |D_i|\le2u_1$, \Cref{lem:strong-boundary} gives
\begin{equation*}\label{eq:component-small-boundary}
 e_G(D_i, D_i^c)-d \ge d+15(M-1)|D_i|\ln d.
\end{equation*}
Plugging this into the last term of \eqref{eq:mixed-master} the total contribution from these sets is at most
\begin{align*}
    &|C \cup D_1| \sum_{k=3}^{2u_1} (ed)^{k-1} (1-p_-)^{\frac{d+15(M-1)k \ln d}{4}}\\
    \leq & n^{3/4} \omega^{1/4} \sum_{k \geq 3} (ed)^{k-1} d^{-\frac{15 \ln 2}{4} k}\\
    \leq & n^{3/4+o(1)},
\end{align*}
where we used that $\ln n \geq \frac{\ln 2}{M-1} d$ and that $\frac{15 \ln 2}{4}>1$.

For $2u_1<|D_i|\le n/2$, first suppose that
\[
 e_G(D_i, D_i^c) \ge \frac{5(M-1)}{\ln 2}|D_i|\ln d.
\]
Using again $\ln n \geq \frac{\ln 2}{M-1}d$, the contribution from this range is at most
\begin{align*}
    &n (1-p_-)^{-d/4}\sum_{k>2u_1}(ed)^k (1-p_-)^{\frac{5(M-1)}{4 \ln 2} k \ln d}\\
    \leq & \left(\frac{\omega}{n} \right)^{-1/4} \sum_{k >2u_1} (ed)^k d^{-\frac{5}{4}k}
 \le\exp\{-\Omega(u_1\ln d)\}.
\end{align*}
Indeed, $\frac{5}{4}>1$ makes the series geometric with exponent $-\Omega(k\ln d)$, while the prefactor has logarithm $O(d)$ and $u_1\gg d$.

Otherwise if $e_G(D_i, D_i^c)<|D_i|\ln^2 d$, then \Cref{lem:small-boundary-sets} (which is applicable because $2u_1 \geq \frac{n}{d^{\ln^2 d}}$) gives at most $\exp\{2|D_i|/\ln d\}$ choices, whereas
\[
 e_G(D_i, D_i^c)-d = \Omega(|D_i|)
\]
since $|D_i|>2u_1 \gg d$. Thus,
\begin{align*}
    &\sum_{k > 2u_1} \exp\left\{ \frac{2k}{\ln d}\right\} (1-p_-)^{\Omega(k)}\\
    \leq & \sum_{k > 2u_1} \exp\{-\Omega(k)\} = \exp\{-\Omega(u_1)\}.
\end{align*}

There is at most one $D_i$ with $|D_i|>n/2$. We now do a similar argument focusing on its complement $D_i^c$. Since $A\subseteq K^c$, we have $u_0<|D_i^c|<n/2$.
If
\[
 e_G(D_i, D_i^c)\ge \frac{8(M-1)}{\ln(2e)}|D_i^c|\ln(en/|D_i^c|),
\]
then,
\[
 \sum_{u_0<k<n/2}\binom nk (1-p_-)^{\left\{\frac{8(M-1)}{\ln(2e)} k\ln(en/k)-d\right\}/4} =\exp(-\Omega(u_0)),
\]
using that $\log n \geq \frac{\log 2}{M-1} d$ and $2 \log 2>1$.
Otherwise if $e_G(D_i, D_i^c)\le |D_i^c|\ln^2d$, then by \Cref{lem:small-boundary-sets} (which is applicable because $u_0 \geq \frac{n}{d^{\ln^2 d}}$) and $e_G(D_i, D_i^c)-d\ge \Omega(|D_i^c|)$. Here again $|D_i^c|>u_0\gg d$ permits the subtraction of $d$ from the linear expansion bound. Then,
\[
\sum_{u_0<k<n/2}\exp\{2k/\ln d\}(1-p_-)^{\Omega(k)}=\exp(-\Omega(u_0)).
\]

Let $\mathcal{S}_{\le n/2}$ be the sum of the terms corresponding to connected sets of size at most $n/2$, and let $\mathcal{S}_{>n/2}$ be the corresponding sum for a possible connected set larger than $n/2$.
The preceding estimates give
\[
 \mathcal{S}_{\le n/2}\le n^{3/4+o(1)},
 \qquad
 \mathcal{S}_{>n/2}\le e^{-\Omega_M(u_0)}.
\]
A pairwise disjoint family contains at most one member larger than $n/2$.
In particular, uniformly for sets $Y$ with $|V\setminus Y|>u_0$, the sum over unordered families of pairwise disjoint connected subsets of $G[Y]$ satisfies
\begin{equation}\label{eq:component-family-sum}
 \sum_{\substack{D^+\\|V\setminus (C \cup D_1)|>u_0}}
 (1-p_-)^{\frac14\sum_{D_i\in D^+}\{e_G(D_i, D_i^c)-d\}}
 \le\exp\{n^{3/4+o(1)}\}
 =\exp\{o(u_0)\}.
\end{equation}
This sum counts only the non-singleton components of $G_p[D]$; it does not use the vertices or internal edges of $C$. Given $C \cup D_1$, an ordered partition $V\setminus (C \cup D_1)=A\sqcup D_0$, and the family $D^+$, the remaining set
$C=(C \cup D_1)\setminus\bigcup_{D_i\in D^+}D_i$ is determined.

Recalling \eqref{eq:mixed-master}, after \eqref{eq:component-family-sum} and after discarding $(1-p_-)^{3dr/8}\le1$, it remains to bound the following sum over ordered partitions $V=(C \cup D_1)\sqcup A\sqcup D_0$:
\begin{equation}\label{eq:ordered-partition-sum}
  \sum_{\substack{V=(C \cup D_1)\sqcup A\sqcup D_0\\
 |A|>u_0}}  (1-p_-)^{e_G((C \cup D_1), (C \cup D_1)^c)/16+\{e_G(A)+e_G(D_0)\}/4}.  
\end{equation}
We split according to whether $s \coloneqq |(C \cup D_1)^c|\le n/2$. Suppose first that this inequality holds. Then $s\ge |A|>u_0$ and $\ln\left(\frac{en}{s}\right)=O(\ln d)$.
If
\begin{equation}\label{eq:small-side-high-boundary}
 e_G((C \cup D_1), (C \cup D_1)^c) \ge \frac{20(M-1)}{\ln 2}s\ln\left(\frac{en}{s}\right),
\end{equation}
then
\[
 \binom ns2^s (1-p_-)^{e_G((C \cup D_1), (C \cup D_1)^c)/16}
 \le \exp\{-\Omega(s\ln(en/s))\}.
\]
If \eqref{eq:small-side-high-boundary} fails, then
\[
 e_G((C \cup D_1), (C \cup D_1)^c)=O(s\ln d)<s\ln^2 d.
\]
By \Cref{lem:small-boundary-sets}, there are $e^{o(s)}$ choices for $(C \cup D_1)^c$, while expansion (\Cref{th:expansion}) gives $e_G((C \cup D_1), (C \cup D_1)^c)= \Omega(s)$.
For each fixed $(C \cup D_1)^c$, \Cref{lem:two-colour-count} gives
\[
 \sum_{A\sqcup D_0}  (1-p_-)^{\{e_G(A)+e_G(D_0)\}/4}\le e^{o(s)},
\]
obtained from \eqref{eq:colouring-weight} with $\theta=1/4$.

Thus, the range $s \leq n/2$ contributes at most
\begin{equation}\label{eq:small-side-sum}
 \sum_{s>u_0}e^{-\Omega(s)}=e^{-\Omega(u_0)}.
\end{equation}
In particular, when the Gallai--Edmonds set $C \cup D_1$ contains more than half of vertices, we count the smaller set $A\cup D_0$ rather than $C \cup D_1$.

We may now assume $s=|(C \cup D_1)^c|>n/2$, so $|C \cup D_1|<n/2$.
Fix $\varepsilon>0$ to be chosen sufficiently small. Suppose first that $|C \cup D_1|\ge\varepsilon n$.
If $e_G((C \cup D_1), (C \cup D_1)^c)\ge \frac{20(M-1) \ln 3}{\ln 2}n$, then
\[
 3^n (1-p_-)^{e_G((C \cup D_1), (C \cup D_1)^c)/16}\le e^{-\Omega(n)}.
\]
Otherwise,
\[
 e_G((C \cup D_1), (C \cup D_1)^c)=O(n)=O(|C \cup D_1|)< |C \cup D_1|\ln^2 d.
\]
By \Cref{lem:small-boundary-sets}, there are $e^{o(n)}$ choices for $(C \cup D_1)^c$, and, for each fixed $(C \cup D_1)^c$, \eqref{eq:colouring-weight} with $\theta=1/4$ (which is applicable since $e_G((C \cup D_1), (C \cup D_1)^c)<\frac{20(M-1) \ln 3}{\ln 2}n\le\frac{40(M-1) \ln 3}{\ln 2}|(C \cup D_1)^c|$) gives
\[
 \sum_{R=A\sqcup D_0}
 (1-p_-)^{\{e_G(A)+e_G(D_0)\}/4}\le e^{o(n)}.
\]
Moreover,
\[
 e_G((C \cup D_1), (C \cup D_1)^c)=\Omega(|C \cup D_1|)\ge \alpha_M\varepsilon n,
\]
for some suitable constant $\alpha_M$.
So, in total this contributes $e^{-\Omega(n)}$.

Finally suppose that $|C \cup D_1|<\varepsilon n$.
Fix $i$ such that $H_i$ is non-bipartite.
Fixing all coordinates other than the $i$-th partitions $V$ into $n/|V(H_i)|$ sets, each of which induces a copy of $H_i$.
At most $|C \cup D_1|$ of these copies meet $C \cup D_1$.
In every remaining copy, the vertices lie in $A\cup D_0$; since $H_i$ is non-bipartite, at least one edge has both endpoints in $A$ or both endpoints in $D_0$.
Thus
\begin{equation}\label{eq:nonbipartite-copies}
 e_G(A)+e_G(D_0)\ge n/M-|C \cup D_1|\ge n/(2M).
\end{equation}
The ordered partitions for which $e_G(A)+e_G(D_0)\ge \frac{5(M-1)\log 3}{\log 2}n$ contribute at most
\[
 3^n (1-p_-)^{\frac{5(M-1)\log 3}{4 \log 2} n}=e^{-\Omega(n)}.
\]
For the remaining ordered partitions with $e_G(A)+e_G(D_0)< \frac{5(M-1)\log 3}{\log 2}n$, repeat the counting argument from \Cref{lem:two-colour-count} for the ordered partition $R=A\sqcup D_0$.
The trivial bound gives $e_G((C \cup D_1), (C \cup D_1)^c)\le d|C \cup D_1|$, so \eqref{eq:size-X} gives only $\varepsilon n+O(n/d)$ exceptional vertices in $X$.
As in the proof of \Cref{lem:two-colour-count}, the random sample $T$, together with $A\mathbin{\triangle}N_{G[R]}(T)$ after possibly interchanging $A$ and $D_0$, determines the ordered partition of $R$. We get
\[
|A \triangle N_{G[R]}(T)| \leq \epsilon n + o(n), \qquad |T|=o(n).
\]
The number of choices for $(T, A \triangle N_{G[R]}(T))$ is at most
\[
 \exp\{\varepsilon\ln\left(\frac{e}{\varepsilon} \right)n+o(n)\}
\]
and the choices of $(C \cup D_1)$ contribute at most $\exp\{\varepsilon\ln(e/\varepsilon)n\}$.
Recalling that by \eqref{eq:nonbipartite-copies} there are at least $(e_G(A)+e_G(D_0))/4 \geq n/(8M)$ edges that are not present in $G_{p_-}$. Choosing $\varepsilon=\varepsilon(M)>0$ sufficiently small ensures $2 \varepsilon \log\left( \frac{e}{\varepsilon}\right) < \frac{\log 2}{8M(M-1)}$. Then, the contribution is at most
\begin{align*}
    &\exp\{\varepsilon\ln\left(\frac{e}{\varepsilon} \right)n+o(n)\} (1-p_-)^{n/(8M)}\\
 \leq &\exp\{2 \varepsilon\ln\left(\frac{e}{\varepsilon} \right)n - \frac{\log 2}{8M(M-1)} +o(n)\}\\
 = & \exp(-\Omega(n)).
\end{align*}
Combining \eqref{eq:small-side-sum} with this contribution from $s>\frac{n}{2}$ proves that \eqref{eq:ordered-partition-sum} is at most $\exp(-\Omega(u_0))+ \exp(-\Omega(n))$.

Multiplying \eqref{eq:component-family-sum} with this bound on \eqref{eq:ordered-partition-sum} gives $o(1)$, showing that whp there is no Gallai--Edmonds decomposition of $G_p$ with $\defi(G_p)\ge2$ and $|A|>u_0$.
\end{proof}

Lemma \ref{lem:large-ge} will follow directly from Lemmas \ref{lem:bipartite-singletons} and \ref{lem:nonbipartite-ge}.

\begin{proof}[Proof of \Cref{lem:large-ge}]
If $G$ is bipartite, every factor-critical bipartite graph is a singleton, so the claim follows from \Cref{lem:bipartite-singletons}.
If $G$ has a non-bipartite base graph, it follows from \Cref{lem:nonbipartite-ge}.
\end{proof}


\section{Discussion}\label{s:discussion}

The deterministic part identifies a natural sufficient condition in terms of the dimension of a product graph to contain a perfect matching that is asymptotically correct. We wonder whether other natural properties or spanning structures are guaranteed by a high dimension of Cartesian product graphs. In particular, we ask the following question for Hamilton cycles.
\begin{question}
    Let $H_i$ for $i \in [t]$ be a sequence of regular graphs of bounded size. Is there a dimension $t^*$ such that for all $t\geq t^*$ the Cartesian product graph $G=\square_{i=1}^t H_i$ has a Hamilton cycle?
\end{question}

In the probabilistic setting we have extended the classical result concerning hitting times of minimum degree one, connectivity, and the existence of a perfect matching to random subgraphs of regular Cartesian product graphs. In particular, this includes a simplified self-contained version of the connectivity result for bond percolation on the hypercube. Let us mention that, independently, Collares, Doolittle, and Erde use a similar approach -- that is, sprinkling with probabilities $p_1$ and $p_2$ -- to show a connectivity result for bond percolation on the permutahedron \cite{CDE24}. There, however, similarly to the approach of \cite{AKS81, BKL92, K23}, one utilises that in $G_{p_1}$, large components are relatively well-spread, that is, typically every vertex in $G$ is quite close (in $G$) to a large component of $G_{p_1}$. In this paper, we neither require nor utilise such a `density' statement, and instead use the fact that the isolated vertices are `sparsely spread'. 

Two features of the proof appear useful beyond this setting.
First, degree summation over the Gallai--Edmonds decomposition turns the matching deficiency into a lower bound on the number of edges that must be absent, expressed through the boundaries of the components of $G_p[D]$. 
Second, the product structure of the graph allows to bound the number of large sets with a small edge-boundary much more efficiently, enabling the use of the union bound, see \Cref{lem:small-boundary-sets}.


Many other random graph models are known to have typically the same hitting times for minimum degree one, connectivity, and the existence of a perfect matching (see, for example, \cite{FP04} and the references therein). It is thus natural to ask what are the minimal requirements on $G$ for this phenomenon to hold. As a step towards this, we propose the following question, considering regular graphs with high-degree.
\begin{question}
    Let $G$ be a $d$-regular graph on $n$ vertices, with $d=\omega(1)$ and $n$ divisible by two. What minimal requirements are needed on $G$, such that in the random graph process on $G$, the hitting times for minimum degree one, connectivity, and the existence of a perfect matching are the same?
\end{question}


\paragraph{Acknowledgements}
The authors would like to thank Joshua Erde, Mihyun Kang, and Michael Krivelevich for their guidance, advice, and fruitful discussions.
Special thanks to an anonymous referee for spotting an error in a previous version.
This research was funded in part by the Austrian Science Fund (FWF) [10.55776/W1230, 10.55776/F1002].
\paragraph{Statement about the use of AI}
The authors acknowledge the use of AI tools for polishing this paper, and note that all of the proof ideas were entirely human-generated.

\begingroup
\hbadness=2000
\bibliographystyle{abbrv}
\bibliography{perc}
\endgroup

\end{document}